\documentclass [a4paper,twoside,10pt]{article}

\usepackage[utf8x]{inputenc}
\usepackage{amsmath,amsfonts,amssymb,amsthm}
\usepackage{vmargin,graphicx}
\usepackage{enumerate}
\usepackage{csquotes}
\usepackage{color}
\usepackage{bm}
\usepackage{soul}
\usepackage[colorlinks,linkcolor=blue,citecolor=blue,urlcolor=blue]{hyperref}

\setpapersize[portrait]{A4}
\setmarginsrb{2cm}{1cm}{2cm}{2.5cm}{1.5cm}{0cm}{0.5cm}{2cm}

\usepackage[english]{babel}
\newcommand{\R}{\ensuremath{\mathbf{R}}}
\newcommand{\B}{\ensuremath{\mathbf{B}}}
\newcommand{\Sp}{\ensuremath{\mathbf{S}}}
\newcommand{\Hp}{\ensuremath{\mathbf{H}}}
\newcommand{\N}{\ensuremath{\mathbf{N}}}

\newtheorem{ethm}{Theorem}[section]
\newtheorem{ecor}[ethm]{Corollary}
\newtheorem{eprop}[ethm]{Proposition}
\newtheorem{elem}[ethm]{Lemma}

\theoremstyle{definition}
\newtheorem{edefi}[ethm]{Definition}

\theoremstyle{remark}
\newtheorem{erem}[ethm]{Remark}

\theoremstyle{theorem}
\newtheorem{thmx}{Theorem}

\newcommand{\proofend}{~$\rhd$}
\newcommand{\proofbegin}{~$\lhd$}

\newenvironment{eproof}
               {\noindent {\emph{\textbf{Proof}}}\\\proofbegin~}
               {\proofend\\}
\newenvironment{eproofof}[1]
               {\noindent {\emph{\textbf{Proof of #1}}}\\\proofbegin~}
               {\proofend\\}

\newcommand{\ABS}[1]{\ensuremath{{\left| #1 \right|}}} 
\newcommand{\PAR}[1]{\ensuremath{{\left(#1\right)}}} 
\newcommand{\SBRA}[1]{\ensuremath{{\left[#1\right]}}} 

\renewcommand{\phi}{\varphi}

\newcommand{\al}{\alpha}

\newcommand{\beq}{\begin{equation}}\newcommand{\eeq}{\end{equation}}
\begin{document}

\title{A conformal geometric point of view on the Caffarelli-Kohn-Nirenberg inequality}

\author{L. Dupaigne, I. Gentil and S. Zugmeyer}

\date{\today}

\maketitle
\begin{abstract}
We are interested in the Caffarelli-Kohn-Nirenberg
inequality (CKN in short), introduced by these authors in 1984. We
explain why the CKN inequality can be viewed as a Sobolev inequality on
a weighted Riemannian manifold. More precisely, we prove that the
CKN inequality can be interpreted in this way on three different and equivalent models, obtained as weighted versions of the
standard Euclidean space, round sphere and hyperbolic space. This result can be viewed as an extension of conformal invariance
to the weighted setting.  Since the spherical CKN model we introduce has
finite measure, the $\Gamma$-calculus introduced by Bakry and \'Emery
provides a way to prove the Sobolev inequalities. This method
allows us to recover the optimality of the region of parameters describing symmetry-breaking of minimizers of the CKN inequality, introduced by Felli and Schneider and proved by
Dolbeault, Esteban and Loss in 2016. Finally, we develop the notion of $n$-conformal invariants, exhibiting a way to extend the notion of scalar curvature to weighted manifolds such as the CKN models.
\end{abstract}

\section{Introduction and main results}
\label{sec-1}
\subsection{The CKN Euclidean space}
In their seminal paper~\cite{ckn}, Caffarelli, Kohn and Nirenberg found the optimal range of real parameters $a,b,p$ for which the following inequality holds true: 
\begin{equation}
  \label{1}
  \PAR{\int_{\R^d}\frac{\vert v\vert ^p}{|x|^{bp}}dx}^{2/p}\leq C_{a,b}\int_{\R^d}\frac{|\nabla v|^2}{|x|^{2a}}dx,\qquad v\in\mathcal C_c^\infty(\R^d\setminus \{0\}).
\end{equation}
Here, $|\cdot|$ is the Euclidean norm in $\R^d$, $d\in\N^*$  and $C_{a,b}$ denotes the optimal constant, depending on $a,b$ and $d$ only. Note that the case $a=b=0$ (and $p=2d/(d-2)$) corresponds to Sobolev's inequality, while the case $b=a+1$ (and $p=2$) is Hardy's inequality, so that~\eqref{1} is sometimes called the Hardy-Sobolev inequality. Note also that the inequality is achieved in the former case, while it is not in the latter. 

Let us consider the measure
\begin{equation}\label{muhat}
  d{\mu_{\bf E}}(x)=|x|^{-bp}dx.
\end{equation}
Then, the left-hand side of~\eqref{1} is simply the $L^p$-norm of $v$ with respect to the measure ${\mu_{\bf E}}$ (squared). In addition, if we consider the metric\footnote{If $(M,g)$ is $d$-dimensional Riemannian manifold whose metric $g$ is represented in a local system of coordinates at a point $x\in M$ by the matrix $G(x)=(g_{ij}(x))_{1\le i,j\le d}$, we use the letter $\mathfrak g$ to denote the bilinear form on the cotangent space of $M$ represented by the inverse matrix $G(x)^{-1}=(\mathfrak g^{ij}(x))_{1\le i,j\le d}$ } on the manifold $M=\R^d\setminus\{0\}$ given by
\begin{equation}\label{ghat}
{\mathfrak g_{\bf E}}^{ij}=|x|^{bp-2a}\delta^{ij},
\end{equation}
then~\eqref{1} takes the simpler form
$$
\PAR{\int \vert v\vert^p d{\mu_{\bf E}}}^{2/p}\leq C_{a,b}\int \vert\nabla_{{\mathfrak g_{\bf E}}}v\vert_{{\mathfrak g_{\bf E}}}^2 \; d{\mu_{\bf E}}.
$$
By a standard scaling argument\footnote{To see this, apply~\eqref{1} to the function $x\mapsto v(\lambda x)$ where $\lambda>0$ and let $\lambda\to 0^+$ and $\lambda\to+\infty$}, the following relation is necessary for the inequality to hold true:
\begin{equation}
\label{crit}
p=\frac{2d}{d-2+2(b-a)}=\frac{d}{a_c-a+b},
\end{equation}
where $a_c=\frac{d-2}{2}$. Through the property of modified inversion symmetry (see Theorem 1.4(ii) in~\cite{florin-wang2001}),  we may always assume that $a<a_c=\frac{d-2}{2}$ since the case $a>a_c$ is dual to it  and the inequality fails to be true if $a=a_c$ (see~\cite{ckn}). For simplicity, we also focus on the case $d\ge4$ and refer to~\cite{DolEstLos14} for the remaining cases $d\in\{1,2,3\}$. Then,~\eqref{1} holds true if and only if 
$$
a\leq b\leq a+1
$$ 
For simplicity, we do not consider the limiting case $b=a+1$ (Hardy's inequality) and we define accordingly the set 
\begin{equation}
\label{468}
\Theta=\{(a,b)\in\R^2,\,a\leq b < a+1,\,a<a_c\}
\end{equation}
so that the CKN inequality~\eqref{1} is valid whenever  $(a,b)\in\Theta$ (see Section~\ref{sec-b}). 

Observe that for $(a,b)\in\Theta$, $p\le\frac{2d}{d-2}$ and so $p$ can be rewritten as the critical Sobolev exponent associated to an {\it intrinsic} dimension $n\in[d,+\infty)$ through the relations
\begin{equation}
\label{dimension}
p=\frac{2n}{n-2},\quad n=\frac{d}{1+a-b}.
\end{equation}
The fact that $n$ is a meaningful number, entering in the classical Bakry-Emery curvature-dimension condition, will become transparent in a moment.  To summarize, one can view inequality~\eqref{1} exactly as Sobolev's inequality stated on the weighted Riemannian manifold\footnote{the words "smooth metric measure space" and  "manifold with density" are also employed in the literature to designate the same object.} that we introduce now. 
\begin{edefi}[The Euclidean CKN space]
\label{eckn} 
The Euclidean CKN space is the triple $(M,{\mathfrak g_{\bf E}},{\mu_{\bf E}})$, where the manifold is $M=\R^d\setminus\{0\}$, the metric\footnote{The given expression of ${\mathfrak g_{\bf E}}$ is just a rewriting of \eqref{ghat}} is ${\mathfrak g_{\bf E}}^{ij}=\vert x\vert^{2(1-\alpha)}\delta^{ij}$\ and where the measure ${\mu_{\bf E}}$ is given by~\eqref{muhat}.
The corresponding Riemannian volume is given by $dV_{{\mathfrak g_{\bf E}}}=\vert x\vert^{d(\alpha-1)}dx$, the weight $W_{\bf E}$, verifying $d{\mu_{\bf E}}=e^{-{W_{\bf E}}}dV_{{\mathfrak g_{\bf E}}}$, is given by ${W_{\bf E}}=-\frac{\alpha(n-d)}2\log\vert x\vert^2$ and the generator\footnote{i.e. the operator such that $-\int u{ L_{\bf E}}v\;d{\mu_{\bf E}}=\int (\nabla_{{\mathfrak g_{\bf E}}}u\cdot_g \nabla_{{\mathfrak g_{\bf E}}}v) d{\mu_{\bf E}}$ for $u,v\in \mathcal C^\infty_c(\R^d\setminus\{0\})$} is given by ${ L_{\bf E}}=\Delta_{{\mathfrak g_{\bf E}}}-\nabla^{{\mathfrak g_{\bf E}}}{W_{\bf E}}\cdot\nabla= \vert x\vert^{2(1-\alpha)}(\Delta -{a}\nabla\log\vert x\vert^2\cdot\nabla)$.
\end{edefi}
For notational convenience, we introduced above the parameter\footnote{The reader may check that $\alpha$ turns out to be the same parameter as the one introduced in~\cite{DolEstLos14} (for different reasons).}:
\begin{equation}
\label{79}
\alpha=1+a-\frac{pb}{2},
\end{equation}
where $(a,b)\in\Theta$ (defined in~\eqref{468})
and $p$ is the critical exponent given by~\eqref{crit}. In other words, returning to the parameters $a,b,d$ (and $a_c=(d-2)/2$),
$$
\alpha=\frac{(a_c-a)(a+1-b)}{a_c-a+b}.
$$
Note that for any $(a,b)\in\Theta$, we have $\alpha\ge 0$, see Section~\ref{sec-b} and Figure~\ref{fig-1} for more information about parameters.

Equivalently, and this is the notation adopted in this paper\footnote{See~\cite{bgl-book} for an introduction to $\Gamma$-calculus.}, one can see the Euclidean CKN space as a Markov triple $(M,{\mu_{\bf E}},{\Gamma_{\bf E}})$, where ${\mu_{\bf E}}$ verifies~\eqref{muhat} and the carré du champ operator is given by
$$
{\Gamma_{\bf E}}(v)=\vert\nabla_{{\mathfrak g_{\bf E}}}v\vert_{{\mathfrak g_{\bf E}}}^2=|x|^{bp-2a}|\nabla v|^2=|x|^{2(1-\alpha)}|\nabla v|^2
$$
Its associated bilinear form is denoted by ${\Gamma_{\bf E}}(u,v)=|x|^{2(1-\alpha)}\nabla u\cdot \nabla v$ for $u,v\in C^\infty_c(\R^d\setminus\{0\})$. 
Inequality~\eqref{1} now reads
$$
\PAR{\int \vert v\vert^p d{\mu_{\bf E}}}^{2/p}\leq C_{a,b}\int {\Gamma_{\bf E}}(v) d{\mu_{\bf E}},\quad v\in\mathcal C_c^\infty(\R^d\setminus \{0\}).
$$

\subsection{Conformal invariance}
As noticed earlier, when $a=b=0$ ($\alpha=1$), we recover the standard Sobolev inequality on the standard Euclidean space. In that case, since the metric of the $d$-dimensional sphere $\Sp^d$ and the metric of the $d$-dimensional hyperbolic space $\Hp^d$ are both conformally equivalent to the Euclidean metric, Sobolev's inequality takes equivalent forms on these three model spaces. More precisely, the Euclidean Sobolev inequality applied to the function $\varphi^{\frac{2-d}2}v$, where $\varphi(x)=\frac{1+\vert x\vert^2}2$ (respectively  $\varphi(x)=\frac{1-\vert x\vert^2}2$) and $v\in \mathcal C^\infty_c(\R^d)$ (resp. $v\in \mathcal C^\infty_c(\B)$) yields 
\begin{equation}\label{confinv}
\PAR{\int \vert v\vert^p dV_{\mathfrak g} }^{2/p}\leq C\left[\int \vert\nabla_{{\mathfrak g}}v\vert_{{\mathfrak g}}^2 \; dV_{\mathfrak g}+\int S_{\mathfrak g} v^2 dV_{\mathfrak g}\right],
\end{equation}
where $\mathfrak g$ is the round metric on the sphere $\Sp^d$ expressed in stereographic cooordinates (resp. the metric of the hyperbolic space in the Poincaré ball model), $dV_{\mathfrak g}$ the associated Riemannian volume, $S_{\mathfrak g}=\frac{d(d-2)}4$ (resp. $S_{\mathfrak g}=-\frac{d(d-2)}4$) and $C=\frac4{d(d-2)}\vert\Sp^d\vert^{-\frac2d}$ the best constant in the standard Euclidean Sobolev inequality. 
By analogy, we can extend the conformal invariance property to the setting of weighted manifolds as described next. 
\subsubsection*{The spherical CKN  and the hyperbolic CKN  spaces}
Recall that the metric and reference measure of the Euclidean CKN space read
$${\mathfrak g_{\bf E}}^{ij}=|x|^{2(1-\alpha)}{\delta^{ij}}\quad\text{and}\quad
d{\mu_{\bf E}}=|x|^{-bp}dx.$$ 

Keeping in mind the expression of the standard stereographic projection, we define next the spherical and hyperbolic CKN spaces as follows.
\begin{edefi}[The spherical and the hyperbolic CKN  spaces]~
\begin{itemize}
\item The spherical CKN space is the triple $(M,{\mathfrak g_{\bf S}}, {\mu_{\bf S}})$, where $M=\R^d\setminus\{0\}$,
$$
{\mathfrak g_{\bf S}}^{ij}=|x|^{2(1-\alpha)}\left(\frac{1+|x|^{2\alpha}}{2}\right)^2{\delta^{ij}}\quad\text{and}\quad 
d{\mu_{\bf S}}={|x|^{-bp}}\left(\frac2{1+|x|^{2\alpha}}\right)^n dx.
$$
Associated objects are given by the following formulae:
\begin{itemize}
\item Riemannian volume: $dV_{\mathfrak g_{\bf S}}=2^d\frac{|x|^{d(\alpha-1)}}{(1+|x|^{2\alpha})^2}dx$,
\item weight: ${W_{\bf S}}=(n-d)\log (1+|x|^{2\alpha})-\frac{\alpha(n-d)}{2}\log |x|^{2},$
\item Carré du champ operator: ${\Gamma_{\bf S}}(v)=|\nabla_{{\mathfrak g_{\bf S}}} v |_{{\mathfrak g_{\bf S}}}^2=|x|^{2(1-\alpha)}\frac{(1+|x|^{2\alpha})^{2}}{4}\vert\nabla v\vert^2$,
\item generator: ${L_{\bf S}}(f)=|x|^{2(1-\alpha)}\frac{(1+|x|^{2\alpha})^2}{4}\SBRA{\Delta f-a\nabla f\cdot\nabla  \log|x|^2-(n-2)\nabla f\cdot\nabla\log(1+|x|^{2\alpha})}.$
\end{itemize}

\item The CKN hyperbolic space is the triple $(\B\setminus\{0\}, {\mathfrak g_{\bf H}},{\mu_{\bf H}})$, where $\B$ is the open unit ball in $\R^d$,
 $$
{\mathfrak g_{\bf H}}^{ij}=|x|^{2(1-\alpha)}\left(
\frac{1-|x|^{2\alpha}}2\right)^2{\delta^{ij}}\quad\text{and}\quad
d{\mu_{\bf H}}=|x|^{-bp}\left(
\frac{2}{1-|x|^{2\alpha}}
\right)^n dx.
$$
Associated objects to this triple are given by the following formulae
\begin{itemize}
\item Riemannian volume:  $d{V_{{\mathfrak g_{\bf H}}}}=2^d\frac{|x|^{d(\alpha-1)}}{(1-|x|^{2\alpha})^{d}}dx$,
\item weight: ${W_{\bf H}}=(n-d)\log (1-|x|^{2\alpha})-\frac{\alpha(n-d)}{2}\log |x|^{2},$
\item Carré du champ operator: ${\Gamma_{\bf H}}(v)=|\nabla_{{\mathfrak g_{\bf H}}} v |_{{\mathfrak g_{\bf S}}}^2=|x|^{2(1-\alpha)}\frac{(1-|x|^{2\alpha})^{2}}{4}\vert\nabla v\vert^2$, 
\item generator: $L_{\bf H}(f)=|x|^{2(1-\alpha)}\frac{(1-|x|^{2\alpha})^2}{4}\SBRA{\Delta f-a\nabla f\cdot\nabla  \log|x|^2-(n-2)\nabla f\cdot\nabla\log(1-|x|^{2\alpha})}.$
\end{itemize}
\end{itemize}
\end{edefi}

\begin{erem}
\begin{itemize}
\item  Note that in the case $\alpha=1$ ({which is achieved in $\Theta$ only when  $a=b=0$, see Lemma~\ref{lem-15}}), the CKN sphere is the standard round sphere (punctured at both of its poles) 
viewed in the stereographic projection chart. Similarily, for $\alpha=1$, the CKN hyperbolic space is the (punctured) hyperbolic space. 
\item Note that, letting
$
\varphi(x)=\frac{1+|x|^{2\alpha}}{2},
$
we have
$
{\Gamma_{\bf S}}=\varphi^2{\Gamma_{\bf E}}\text{ and }{\mu_{\bf S}}=\varphi^{-n}{\mu_{\bf E}}.
$
We shall say that the CKN Euclidean and spherical spaces belong to the same $n$-conformal class ($n$ not necessarily being equal to the topological dimension). Similarly, with $
\psi(x)=\frac{1-|x|^{2\alpha}}{2},
$
we have
$
{\Gamma_{\bf H}}=\psi^2{\Gamma_{\bf E}}\text{ and }{\mu_{\bf H}}=\psi^{-n}{\mu_{\bf E}},
$
so that the hyperbolic CKN space also belongs to the same $n$-conformal class.
\item When $(a,b)\in\Theta$, ${\mu_{\bf S}}$ has finite mass (see Remark \ref{mufini} in Section~\ref{sec-b}). In this Section, we prefer not to normalize the measure ${\mu_{\bf S}}$, a choice which makes the conformal invariance of Sobolev's inequality more transparent.  
\end{itemize}
\end{erem}

With these definitions at hand, we prove 
\begin{ethm}[Conformal invariance of the three model spaces]
\label{confinvckn}
Let $C>0$ be an arbitrary constant. The three following Sobolev inequalities associated to each CKN model are equivalent:
\begin{align}
  (i) \quad&   \forall v\in C^\infty_c(\R^d\setminus\{0\}),
  &&\PAR{\int \vert v\vert^p d{\mu_{\bf E}} }^{2/p}\leq C\int {\Gamma_{\bf E}}(v) d{\mu_{\bf E}}, \label{179a}\\
  (ii) \quad&  \forall v\in C^\infty_c(\R^d\setminus\{0\}),
  &&\PAR{\int \vert v\vert^p d{\mu_{\bf S}} }^{2/p}\leq C\left(\int {\Gamma_{\bf S}}(v) d{\mu_{\bf S}}+\frac{n(n-2)}{4}\alpha^2\int v^2 d{\mu_{\bf S}}\right),\label{179b}\\
  (iii) \quad& \forall v\in C^\infty_c(\B\setminus\{0\}),
  &&\PAR{\int \vert v\vert^p d{\mu_{\bf H}} }^{2/p}\leq C\left(\int {\Gamma_{\bf H}}(v) d{\mu_{\bf H}}-\frac{n(n-2)}{4}\alpha^2\int v^2 d{\mu_{\bf H}}\right). \label{179c}
\end{align}
\end{ethm}
 
\begin{erem}
  \begin{itemize}
  \item Inequality \eqref{179a} is valid for some constant \(C=C_{a,b}\) if and only if $(a,b)\in\Theta$ as proved in \cite{ckn}. Hence, so are \eqref{179b} and \eqref{179c}.
  \item As we shall see, the value of the optimal constant $C$ is known only in a restricted range of parameters, see Theorem \ref{thm-1} below.
  \item Since its set of test functions is smaller, inequality \eqref{179c} need not be optimal even though \eqref{179a} and \eqref{179b} are. For example, in the absence of weights, $C=\frac4{d(d-2)\vert\Sp^d\vert^{2/d}}$ is the optimal constant in \eqref{179a} and \eqref{179b} and extremals exist (and are classified), see Theorem \ref{thm-1} below. In contrast, inequality \eqref{179c} holds with the same constant $C=\frac4{d(d-2)\vert\Sp^d\vert^{2/d}}$ but when $n=d=3$ (and again $\alpha=1$), the constant $-\frac{n(n-2)}{4}$ can be improved to $-\frac{(n-1)^2}4$, see~\cite{benguria}.
Using this fact and the proof of Theorem \ref{confinvckn}, it follows that the standard Sobolev inequality in $\R^3$ improves to
$$
\PAR{\int_{\R^3} \vert v\vert^{6} dx }^{1/3}\leq \frac13 \left(\frac2{\pi}\right)^{4/3}\left(\int_{\R^3} \vert\nabla v\vert^2 dx - \int_{\R^3}\frac{v^2}{(1-\vert x\vert^2)^2}dx\right),
$$
when restricted to functions $v\in H^1_0(\B)$. When $d\ge4$ (and $\alpha=1$), inequality \eqref{179c} is again optimal, but contrary to~\eqref{179a} and \eqref{179b}, the inequality is never attained\footnote{If it were, then \eqref{179a} would also be attained by a compactly supported function.}. We have not yet investigated the optimality of~\eqref{179c} in the general  case.
  \end{itemize}
 \end{erem}

\subsection{Curvature-dimension conditions on the spherical CKN space}

We just saw that the CKN inequality takes the forms~\eqref{179a},  \eqref{179b}, \eqref{179c} on the three CKN spaces. 
 But what is the value of the best constant $C$ ?  To answer this question, let us first recall the following classical definition and result: a smooth weighted manifold $(M,\mathfrak g,\mu)$ is said to satisfy the $CD(\rho,n)$ condition if for every $f\in \mathcal C^\infty(M)$,
\begin{equation*}
\Gamma_2(f)\ge \rho\Gamma(f)+\frac1n(Lf)^2,
\end{equation*}
where $(\rho,n)\in\R\times(\R\cup \{+\infty\})$, $d\mu=e^{-W}dV_g$ for some $W\in \mathcal C^\infty(M)$,  $\Gamma(f,h)=\langle\nabla^{\mathfrak g}f,\nabla^{\mathfrak g}h\rangle_{\mathfrak g}$, $\Gamma(f)=\Gamma(f,f)=\vert\nabla^{\mathfrak g}f\vert_{\mathfrak g}^2$, $Lf=\Delta_{\mathfrak g}f-\Gamma(W,f)$ and $\Gamma_2(f) = \frac12  L(\Gamma(f))-\Gamma(f, Lf)$, for smooth functions $f,g$ on $M$.
The following theorem, generalizing the earlier work \cite{BidVer91}, holds true:

\begin{thmx}{\rm (\cite{BakLed96} and \cite[Thm.~6.8.3]{bgl-book})}
\label{th:B}
 Let $(M,\mathfrak g,\mu)$ be a smooth weighted manifold satisfying the $CD(\rho,n)$ condition with $\rho>0$, $n>d=dim(M)$ ($n>2$). Assume in addition that the associated operator $L$ is essentially self-adjoint in $L^2(\mu)$.
Let $p=\frac{2n}{n-2}$ and normalize the measure $\mu$ so that $\mu(M)=1$. Then,
\begin{equation}
\label{bl96}
\PAR{\int \vert v\vert^p d\mu }^{2/p}\leq  \frac{4}{n(n-2)}\frac{n-1}\rho\int \Gamma(v) d\mu+\int v^2 d\mu, \quad v\in\mathcal C_c^\infty(M).
\end{equation}	
\end{thmx}
Let us remark that this theorem can also be stated in the more general context of full Markov triples as proposed in~\cite{bgl-book} and also on  metric measure spaces as proved in~\cite{profeta}.
Thanks to Theorem~\ref{th:B}, it suffices to determine whether the CKN sphere is a smooth weighted manifold  satisfying the $CD(\rho,n)$ condition and that the associated operator is essentially self-adjoint, in order to obtain an explicit value (which turns out to be optimal in our case) for the constant $C$ in~\eqref{179b}. This is what we do next.
\begin{eprop}[Curvature-dimension condition for the spherical CKN space]
\label{prop-100}
Let $(a,b)\in\Theta$ and 
\begin{equation}
\label{360}
\rho=\alpha^2(n-1).
\end{equation}
Then, the spherical CKN space  satisfies the curvature-dimension condition $CD(\rho, n)$ if and only if
\begin{equation}\label{bdgz}
\alpha^2\le\frac{d-2}{n-2}.
\end{equation}
\end{eprop}
One implication of this proposition has been proved in~\cite[Thm.~3.9]{Ketterer}. The proof proposed  here is different and based on tensors, which is a useful method to prove the equivalence between the two conditions.  
Note that $(\R^d\setminus\{0\},{\mathfrak g_{\bf S}},{\mu_{\bf S}})$ is a smooth weighted manifold and that its operator is essentially self-adjoint  if and only if $n\ge 3$,
see~\cite[Thm.~3.12]{Ketterer}, whence Sobolev's inequality~\eqref{179b} holds under the condition~\eqref{bdgz}. In fact, more can be said. Revisiting the proofs of Theorem~\ref{th:B} given in~\cite{dgz, bgl-book}, we find that if the following {\it weaker} integrated form of the curvature-dimension condition 
\begin{equation}\label{cdrhoi}
\int \left(\Gamma^{\bf S}_2(f)- \rho\Gamma_{\bf S}(f)-\frac1n(L_{\bf S}f)^2\right)f^{1-n}d\mu_{\bf S}\ge 0
\end{equation}
holds for functions $f\in \mathcal C^\infty(\R^d\setminus\{0\})$ such that $\inf f>0$ and $\sup f<+\infty$, then the sharp Sobolev inequality is valid on the CKN sphere. More precisely, letting $H^1_0(\mu_{\bf S})$ denote the closure of $C^\infty_c(\R^d\setminus\{0\})$ with respect to the norm 
$$
\Vert f\Vert_{H^1_0({\mu_{\bf S}})}^2 = \int ({\Gamma_{\bf S}}(f)+f^2)\;d{\mu_{\bf S}}.
$$


\begin{ethm}[Sobolev inequality for the spherical CKN space]
\label{thm-1}
Let $(a,b)\in\Theta$. Whenever 
\begin{equation}\label{fs}
0<\alpha< 1,
\end{equation}
the following optimal Sobolev inequality holds
	\begin{equation}
	\label{355}
	\PAR{\int v^pd\mu}^{2/p}\leq \frac{4}{n(n-2)\alpha^2}\int{\Gamma_{\bf S}}(v)d\mu+\int v^2d\mu,
	\end{equation}
	for any $v\in H^1_0({\mu_{\bf S}})$, where $\mu=\frac1Z{\mu_{\bf S}}$ {and $Z$ is a normalization constant\footnote{$Z={\mu_{\bf S}}(\R^d\setminus\{0\})=\frac2\alpha\vert \Sp^{d-1}\vert\int_0^{+\infty}(\cosh t)^{-n}dt$} such that $\mu$ is a probability measure}. 
	That is, inequality~\eqref{179b} is valid with optimal constant  
	$$
	C=\frac{4}{n(n-2)\alpha^2Z^{\frac2n}}.
	$$ 
In addition, equality holds in~\eqref{355} if and only if 
$$
v(x)=(\lambda+\gamma\tanh(\alpha s))^{-\frac{n-2}{2}},\quad s=\log\vert x\vert,
$$
where $\lambda,\gamma$ are arbitrary constants such that $\lambda>\vert\gamma\vert$.
In particular optimal functions for both inequalities~\eqref{1} and~\eqref{355} are radial.
\end{ethm}

\begin{erem}

\begin{itemize}
\item Note that $\varphi_1(x)=\tanh(\alpha s)$ is a radial eigenfunction of ${L_{\bf S}}$ associated to the eigenvalue $\lambda=\alpha^2 n$. So, except for the round sphere (corresponding to the case \(\alpha=1\)), the extremals of Sobolev's inequality are obtained as a linear combination of radial extremals of Poincaré's inequality  \eqref{890} (raised to the power $-\frac{n-2}2$), provided this combination is bounded below by a positive constant.
\item  As we shall prove in Lemma~\ref{lem-15}, condition~\eqref{fs} is equivalent to \(\alpha^2\leq \frac{d-1}{n-1}\). In the limiting case $\alpha^2=\frac{d-1}{n-1}$, extremals of Sobolev's inequality are radial, while extremals of Poincaré's inequality need not be, see Proposition \ref{poincare} below.
\item The extremals of Sobolev's inequality on the round sphere (i.e. the limiting case $\alpha=1$) were discovered by T. Aubin, see e.g. Theorem 5.1 in~\cite{Heb99}. They are more often written as constant multiples of
$$
v=(\beta-\cos(r))^{-\frac{d-2}{2}},
$$
where $\beta>1$ and $r$ is the geodesic distance to an arbitrary point $\omega_0\in\Sp^d$. With our notations, they take the form 
$$
v=(\lambda+\gamma\varphi_{1,d})^{-\frac{d-2}{2}},
$$
where $\varphi_{1,d}$ is any eigenfunction of $-{L_{\bf S}}=-\Delta_{\Sp^d}$ associated to the first nonzero eigenvalue $\lambda_1=d$ and $\lambda,\gamma$ are arbitrary constants such that $\lambda>\vert\gamma\vert\Vert\varphi_{1,d}\Vert_\infty$.

Our notation puts forward the connection between the extremals of Sobolev's inequality and the extremals of Poincaré's inequality in Proposition~\ref{poincare}:  the former are obtained as a linear combination of the latter (raised to the power $-\frac{d-2}2$), provided this combination is bounded below by a positive constant.

\item Hardy's inequality (i.e. the case $\alpha=0$ in \eqref{179a}) is optimal for the constant 
$$
C=\lim_{\alpha\to0^+}\frac{4}{n(n-2)\alpha^2Z^{\frac2n}}=\left(\frac2{d-2-2a}\right)^2,
$$
but equality is never achieved.
\item It follows from Theorem \ref{confinvckn} and Theorem \ref{thm-1} that for $\alpha\in(0,1)$, extremal functions for Sobolev's inequality on the Euclidean CKN space take the form
$$
v(x)= \left(\frac{1+|x|^{2\alpha}}{2}\right)^{-\frac{n-2}2},
$$
up to normalization and dilation, providing thereby an alternative proof of the main result in \cite{DolEstLos14}.
\end{itemize}

\end{erem}

\begin{figure}[!h]
\begin{center}
\includegraphics[width=10cm]{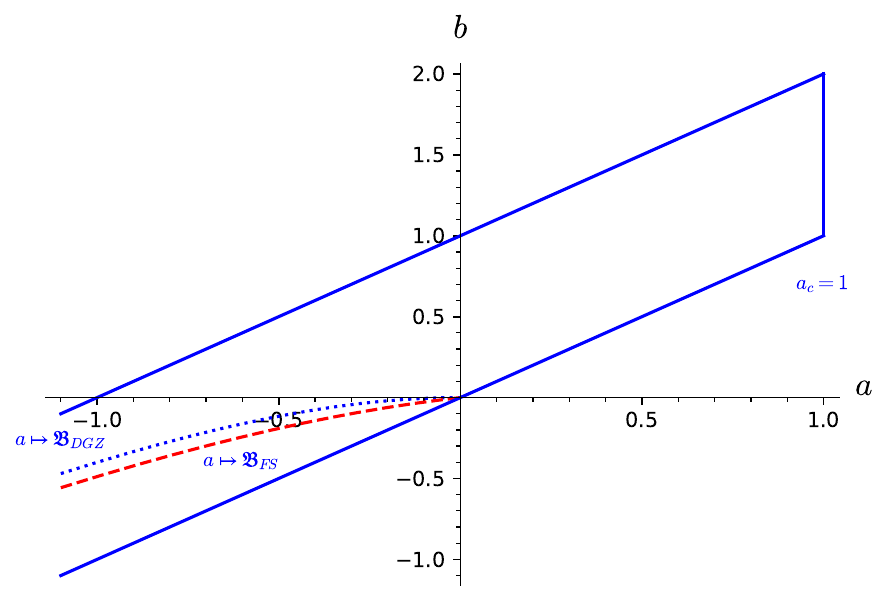}
\caption{$b_{FS}(a)$  (dashed style in red) and the curve $\mathfrak B_{DGZ}=0$ (dotted style in blue) with  $d=4$.}
\label{fig-1}
\end{center}
\end{figure}

Condition~\eqref{fs} is strictly weaker than condition~\eqref{bdgz}. It turns out to be equivalent to 
$$
\alpha^2\le\frac{d-1}{n-1}.
$$
More precisely, consider the following Felli-Schneider region 
\begin{equation}
\label{865}
\Theta_{FS}=\{(a,b)\in\Theta,\,\, b\geq b_{FS}(a)\,\,{\rm if}\,\, b\leq0\},
\quad\text{
where}\quad
b_{FS}(a)=\frac{d(a_c-a)}{2\sqrt{(a_c-a)^2+d-1}}-(a_c-a)
\end{equation}
Let as well \(\Theta_{DGZ}\subset\Theta\) be the domain where Proposition~\ref{prop-100} is valid (condition~\eqref{bdgz}):
$$
\Theta_{DGZ}=\left\{(a,b)\in \Theta,\, \mathfrak B_{DGZ}\geq0\right\},
$$
 where 
\begin{equation}
\label{951}
\mathfrak B_{DGZ}=(d-2)-(n-2)\alpha^2.
\end{equation}
Then, we prove that the region \(\Theta_{FS}\) corresponds exactly to the domain where Theorem~\ref{thm-1} is valid:
\begin{elem}[Comparison of the two regions]
\label{lem-15}
$$
\Theta_{FS}=\left\{(a,b)\in \Theta,\, \mathfrak B_{DGZ}+\frac{n-d}{n-1}\geq0\right\}=\left\{(a,b)\in \Theta,\, \alpha^2\leq \frac{d-1}{n-1}\right\}=\left\{(a,b)\in \Theta,\,\alpha\in(0,1]\right\}.
 $$
Hence, $\Theta_{DGZ}\varsubsetneq \Theta_{FS}$. Moreover, for any $(a,b)\in\Theta$, $\alpha=1$ if and only if $a=b=0$ (that is the standard Euclidean case).
\end{elem}
The two regions are represented in Figure~\ref{fig-1} with $d=4$.
The condition $\theta\in\Theta_{FS}$ is known in the litterature as the Felli-Schneider condition. 
Felli and Schneider~\cite{felli-schneider}, building on the work of Catrina and Wang~\cite{florin-wang2001} initially proved that extremal functions for the optimal CKN inequality~\eqref{1} cannot be radial whenever~\eqref{fs} fails. Conversely, in their work~\cite{DolEstLos14}, Dolbeault, Esteban and Loss computed the optimal constant in~\eqref{1} and proved that extremals for the optimal CKN inequality~\eqref{1} are radial and explicit whenever~\eqref{fs} holds. Combining Theorem~\ref{thm-1} with Theorem~\ref{confinvckn} gives an immediate alternative proof of these latter facts. Our point of view may further clarify why the Felli-Schneider condition is optimal. Indeed, it is well-known that a tight Sobolev inequality implies a Poincaré inequality: precisely, applying~\eqref{355} with $v=1+\epsilon f$ and letting $\epsilon\to0$ leads to
\begin{equation}
\label{890}
\int f^2d\mu-\PAR{\int fd\mu}^2\leq \frac{1}{n\alpha^2}\int {\Gamma_{\bf S}}(f)d\mu.
\end{equation}
We prove the following.
\begin{eprop}[Poincaré inequality for the spherical CKN space]\label{poincare}
Let $(a,b)\in\Theta$. The Poincaré inequality~\eqref{890} holds with optimal constant $C=\frac{1}{n\alpha^2}$ if and only if the Felli-Schneider condition~\eqref{fs} holds.  In addition, equality holds in \eqref{890} if and only if $f$ is an eigenfunction associated to the first nonzero eigenvalue of ${L_{\bf S}}$,  
\begin{itemize}
\item if $0<\alpha<\frac{d-1}{n-1}$, 
$$
f(x)=\lambda+\gamma\tanh(\alpha s), \quad\text{where $s=\log{\vert x\vert}$,}
$$ 
for some constants $(\lambda,\gamma)\in\R^2$.
\item otherwise, if $\alpha=\frac{d-1}{n-1}$,
$$
f(x)=\lambda+\gamma\tanh(\alpha s)+\nu\frac{\varphi_{1,d-1}(\omega)}{\cosh(\alpha s)}, \quad\text{where $s=\log{\vert x\vert}$ and $\omega=\frac{x}{\vert x\vert}$,}
$$ 
with  $(\lambda,\gamma,\nu)\in\R^3$ and some eigenfunction $\varphi_{1,d-1}$ associated to the first nonzero eigenvalue $\lambda_1=d-1$ of $-\Delta_{\Sp^{d-1}}$.
\end{itemize}
\end{eprop}

So, the Felli-Schneider condition cannot be improved in the statement of Theorem~\ref{thm-1} and it is in fact equivalent to both Sobolev's and Poincaré's inequality with the given optimal constants on the spherical CKN space. In addition, Poincaré's inequality with constant $C=\frac{n-1}{\rho n}$ is in fact equivalent to the following integrated $CD(\rho,n)$ condition, where $\rho,n>0$:
$$
\int \left(\Gamma_2^{\bf S}(f)- \rho\Gamma_{\bf S}(f)-\frac1n(L_{\bf S}f)^2\right)d\mu\ge 0,
$$
see Proposition 4.8.3, Theorem 4.8.4 and their proofs in~\cite{bgl-book}. Hence, the Felli-Schneider condition~\eqref{fs} can be interpreted as a curvature-dimension condition in integral form. 

Let us also point out that if Sobolev's inequality~\eqref{179b} holds on a $d$-dimensional smooth manifold $(M,\mathfrak g)$ {\it without weight} with optimal constant $C=\frac{4}{d(d-2)\vert \Sp^{d}\vert^{2/d}}$, then $(M,\mathfrak g)$ must be isometric to the round sphere, as very recently demonstrated in~\cite{nobili}. 


\subsection{The \boldmath{$n$}-conformal invariants}
In this last introductory paragraph, we expand on the conformal invariance of Sobolev's inequality in the setting of weighted manifolds and provide a deeper reason for why the three CKN model spaces satisfy equivalent conformal forms of the Sobolev inequality. 

For the inequality~\eqref{confinv} without weights, it turns out that $S_{\mathfrak g}=  \frac{d-2}{4(d-1)}sc_{\mathfrak g} $ is a constant multiple of the scalar curvature of $\mathfrak g$ (see Propositions~3.6.20,~3.6.21 and~6.2.2, as well as the second displayed formula on p.~63 in~\cite{Heb97} or~\cite[Sec.~6.9.2]{bgl-book} for proofs of this classical result). In other words, inequality \eqref{confinv} is valid on the whole conformal class of the round sphere, including the Euclidean (where $S_{\mathfrak g}=0$) and hyperbolic (where $S_{\mathfrak g}=  -\frac{d-2}{4(d-1)}sc_{\mathfrak g} $ ) spaces. 

For weighted manifolds, the notion of scalar curvature can be generalized as follows. As proposed in~\cite[Sec.~6.9]{bgl-book} (see also~\cite{chang} and~\cite{case} for earlier perspectives\footnote{which correspond to the special case $\gamma=-2$ in Proposition~\ref{prop1} below.}), given a $d$-dimensional ($d\geq 2$) weighted Riemannian manifold $(M,\mathfrak g,\mu)$ with reference measure 

$$d\mu=e^{-W}dV_{\mathfrak g},$$ 

where $W:M\to\R$ is a given weight and  $dV_{\mathfrak g}$ the Riemannian volume, let  
$$
\Gamma(f)=|\nabla_{\mathfrak g} f |_{\mathfrak g}^2
$$ 
denote the associated carré du champ operator, so that $(M,\mu,\Gamma)$ is a Markov triple. 
\begin{edefi}
Take a real number $n\in[d,+\infty]$, which is not necessarily an integer.  The $n$-conformal class of the triple $(M,\mu,\Gamma)$ is the set of all Markov triples $(M, c^{-n}\mu, c^2 \Gamma)$, where $c:M\mapsto (0,\infty)$ is any smooth and positive function.  An $n$-conformal invariant is a map $S$ defined on the $n$-conformal class of  $(M,\mu,\Gamma)$ with values in the set of functions over $M$,
such that for any positive smooth function $c=e^\tau$, 
\begin{equation}
\label{178}
S(c^{-n}\mu,c^2\Gamma)=c^2\SBRA{S(\mu,\Gamma)+\frac{n-2}{2}\left(L\tau-\frac{n-2}{2}\Gamma(\tau)\right)},
\end{equation}
where $$L=\Delta_{\mathfrak g}-\Gamma(W,\cdot).$$
\end{edefi}
It is important to notice that the operator $L$ is uniquely determined by the carré du champ operator $\Gamma$ and the measure $\mu$ only. This is indeed the case since the operator $\Delta_{\mathfrak  g}$ depends on the metric $\mathfrak g$ (which itself is uniquely determined by $\Gamma$) and since $W$ is related to the measure $\mu$ and the metric $\mathfrak g$ through the Riemannian measure $dV_{\mathfrak g}$. Also observe that setting $u=c^{-\frac{n-2}{2}}$, $s=S(\mu,\Gamma)$, $\tilde s=S(c^{-n}\mu,c^2\Gamma)$, then \eqref{178} can be reformulated as the following Yamabe-type equation:
$$
-Lu + su = \tilde s u^{\frac{n+2}{n-2}}\quad\text{in $M$.}
$$
Note that the case where $n=d$, $L=\Delta_{\mathfrak g}$, $s=\frac{d-2}{4(d-1)}sc_{\mathfrak g}$ and $\tilde s$ constant, is the standard Yamabe equation.
\medskip

By a rather direct computation, see~\cite[Prop.~6.9.2]{bgl-book}, whenever $S=S(\mu,\Gamma)$ is an $n$-conformal invariant, the  Sobolev inequality 
\begin{equation}
\label{179d}
\PAR{\int \vert v\vert^p d\mu }^{2/p}\leq C\left(\int \Gamma(v) d\mu+\int S v^2 d\mu\right),
\end{equation}
(with given constant $C>0$ and $p=\frac{2n}{n-2}$)
is invariant in the $n$-conformal class of the triple $(M,\mu,\Gamma)$. In other words, if the Sobolev inequality~\eqref{179d} holds for some constant $C$, then it also holds with the same constant $C$ for all triples $(M, c^{-n}\mu, c^2 \Gamma)$ where $c$ is any smooth and positive function. 

\medskip

Rephrasing what we said earlier, in the absence of weight, $S_{\mathfrak g} = \frac{d-2}{4(d-1)}sc_{\mathfrak g}$ is an example of a $d$-conformal invariant (where in this case $n=d$). The case of weighted Riemannian manifolds is a little bit more complicated and contains interesting examples.  Let us recall Proposition~6.9.6 of~\cite{bgl-book} (we will also provide a proof since the one in~\cite{bgl-book} contains some mistakes as well as the statement). 

\begin{eprop}[$n$-conformal invariant in a weighted manifold]
\label{prop1}	
Let  $\gamma\in\R$ and $n>d$. Then, 
\begin{equation}
\label{181}
S_\gamma(\mu,\Gamma)=\theta_n(\gamma)[sc_{\mathfrak g}-\gamma\Delta_{\mathfrak g}W+\beta_n(\gamma)\Gamma(W)]
\end{equation}
is an $n$-conformal invariant if
$$
\beta_n(\gamma)=\frac{\gamma(n-2d+2)-2(d-1)}{2(n-d)}\quad\text{and}\quad \theta_n(\gamma)=\frac{n-2}{4(d-1)-2\gamma(n-d)}.
$$
\end{eprop}
This being recalled, a natural question arises in the context of the Euclidean CKN space we introduced in Definition~\ref{eckn}: does there exist a (unique) real number $\gamma_0\in\R$ such that this space satisfies 
$$S_{\gamma_0}({\mu_{\bf E}},{\Gamma_{\bf E}})=0\quad ?$$   
This is indeed the case, as we are about to see. By Theorem~\ref{confinvckn}, without any further computation, we deduce that for the same value of the parameter $\gamma=\gamma_0$,  $S_{\gamma_0}({\mu_{\bf S}},{\Gamma_{\bf S}})=\frac{n(n-2)}{4}\alpha^2>0$ is constant for the CKN sphere and  $S_{\gamma_0}({\mu_{\bf H}},{\Gamma_{\bf H}})=-\frac{n(n-2)}{4}\alpha^2$ for the CKN hyperbolic space. 
\begin{eprop}
	\label{prop-5}
 Let $n>d$ and 
\begin{equation}\label{gamma}
\gamma_0=-\frac{2(d-1)}{\alpha^2(n-d)(n-2)}\mathfrak B_{DGZ},
\end{equation}
where $\mathfrak B_{DGZ}$ has been defined in~\eqref{951}.

Then, $S_{\gamma_0}({\mu_{\bf E}},{\Gamma_{\bf E}})=0$, $
	S_{\gamma_0}({\mu_{\bf S}},{\Gamma_{\bf S}})= \frac{n(n-2)}{4}\alpha^2
	$
and 
$
S_{\gamma_0}({\mu_{\bf H}},{\Gamma_{\bf H}})= - \frac{n(n-2)}{4}\alpha^2.
$
\end{eprop}
As an immediate collorary of the CKN inequality~\eqref{1} and the above lemma, we recover the validity of Sobolev's inequality on our three model spaces, stated in Theorem~\ref{confinvckn} above. 

\begin{erem}
In a forthcoming report, we will further explain how a weighted version of Otto's calculus can be introduced in order to prove a wider class of optimal CKN inequalities, by working directly on the Euclidean CKN space, rather than the CKN sphere. 
\end{erem}

{The rest of the paper is organized as follow. 
In Section 2 below, we prove the conformal invariance of Sobolev's inequality in the CKN spaces (Theorem~\ref{confinvckn}). Section 3 is dedicated to the characterization of the region of parameter $\Theta_{DGZ}$ (resp. $\Theta_{FS}$) for which the classical curvature-dimension condition (resp. the integrated form \eqref{cdrhoi}) holds, from which Sobolev's inequality follows (Proposition~\ref{prop-100} and Theorem~\ref{thm-1}). In Section 4, we prove all results pertaining to $n$-conformal invariance for general weighted manifolds (Propositions~\ref{prop1} and~\ref{prop-5}). At last, an appendix contains lists of known formulas and constants, proofs of the numerology relating them as well as rigorous justification of the integrations by parts implicitly used in the proof of Sobolev's inequality.}

\section{Conformal invariance of Sobolev type inequalities for CKN models}
\label{sec-4}

 {\noindent {\emph{\textbf{Proof of Theorem~\ref{confinvckn}}}}\\\proofbegin~}
As in the unweighted case, the proof reduces to a simple change of unknown, once the proper notion of conformal invariance has been introduced. First we prove that the Sobolev inequality in the CKN Euclidean space is equivalent to the Sobolev inequality in the spherical CKN space.

Recall that for
\begin{equation}\label{eq:phi}
\varphi(x)=\frac{1+|x|^{2\alpha}}{2},\quad x\in\R^d,
\end{equation}
and that 
$$
{\Gamma_{\bf S}}=\varphi^2{\Gamma_{\bf E}}\quad\text{and}\quad{\mu_{\bf S}}=\varphi^{-n}{\mu_{\bf E}}. 
$$
Apply~\eqref{1} to the function $f=\varphi^{\frac{2-n}{2}}g$.  On the one hand, we have
$$
\int f^pd{\mu_{\bf E}}=\int f^{\frac{2n}{n-2}}d{\mu_{\bf E}}=\int g^pd{\mu_{\bf S}}. 
$$
On the other hand, letting $V=\log \varphi$,
\begin{align*}
{\Gamma_{\bf E}}(f)={\Gamma_{\bf E}}(\varphi^{\frac{2-n}{2}}g)&=\varphi^{2-n}{\Gamma_{\bf E}}(g)+2\varphi^{\frac{2-n}2}g{\Gamma_{\bf E}}(\varphi^{\frac{2-n}2},g)+{\Gamma_{\bf E}}(\varphi^{\frac{2-n}2})g^2\\
&=\varphi^{2-n}\PAR{{\Gamma_{\bf E}}(g)-\frac{n-2}{2}{\Gamma_{\bf E}}(g^2,V)+\frac{(n-2)^2}{4}{\Gamma_{\bf E}}(V)g^2}.
\end{align*}
An integration by parts with respect to \({\mu_{\bf E}}\) yields
\[
\int {\Gamma_{\bf E}}(g^2,V)\varphi^{2-n}d{\mu_{\bf E}} + \int {\Gamma_{\bf E}}(\varphi^{2-n},V)g^2d{\mu_{\bf E}} = -\int { L_{\bf E}}V   g^2\varphi^{2-n}d{\mu_{\bf E}},
\]
so that we get 
\begin{align*}
\int{\Gamma_{\bf E}}(f)d{\mu_{\bf E}}&=\int{\Gamma_{\bf S}}(g)d{\mu_{\bf S}}-\frac{n-2}{2}\int {\Gamma_{\bf E}}(g^2,V)\varphi^{2-n}d{\mu_{\bf E}}+\frac{(n-2)^2}{4}\int{\Gamma_{\bf E}}(V)g^2\varphi^2 d{\mu_{\bf S}}\\
&=\int{\Gamma_{\bf S}}(g)d{\mu_{\bf S}}+\frac{n-2}{2}\int\PAR{{ L_{\bf E}}(V)-\frac{n-2}{2}{\Gamma_{\bf E}}(V)} g^2\varphi^2 d{\mu_{\bf S}}, 
\end{align*}
The CKN inequality~\eqref{1} becomes 
$$
\PAR{\int g^pd{\mu_{\bf S}}}^{2/p}\leq { }C\left(\int {\Gamma_{\bf S}}(g) d{\mu_{\bf S}}+\frac{(n-2)}{2}
\int \PAR{{ L_{\bf E}}(V)-\frac{n-2}{2}{\Gamma_{\bf E}}(V)}\varphi^2g^2d{\mu_{\bf S}}\right),
$$
and it is enough to compute the quantity 
$$
{L_{\bf E}}(V)-\frac{n-2}{2}{\Gamma_{\bf E}}(V).
$$
To that end, recall that \({ L_{\bf E}}\) is given by
\begin{equation}
\label{47}
L_{\bf E}=|x|^{2(1-\alpha)}\SBRA{\Delta -a\nabla \log |x|^2\cdot\nabla }.
\end{equation}
Since $V=\log \varphi$, we have 
$$
{L_{\bf E}}(V)-\frac{n-2}{2}{\Gamma_{\bf E}}(V)=\frac{{ L_{\bf E}}(\varphi)}{\varphi}-\frac{{\Gamma_{\bf E}}(\varphi)}{\varphi^2}-\frac{n-2}{2}\frac{{\Gamma_{\bf E}}(\varphi)}{\varphi^2}=\frac{{ L_{\bf E}}(\varphi)}{\varphi }-\frac{n}{2}\frac{{\Gamma_{\bf E}}(\varphi)}{\varphi^2}.
$$
Recalling the definition of $\varphi$ given in~\eqref{eq:phi}, we find that
$$
\frac{{\Gamma_{\bf E}}(\varphi)}{\varphi^2}=\frac{\alpha^2|x|^{2\alpha}}{\varphi^2},
$$
and
$$
\frac{L_{\bf E}(\varphi)}{\varphi}=\frac{|x|^{2(1-\alpha)}\SBRA{\Delta \varphi-a\nabla \log |x|^2\cdot\nabla \varphi}}{\varphi}=\alpha\frac{d+2(\alpha-1)-2a}{\varphi}=\frac{\alpha^2n}{\varphi}.
$$
So finally,  
$$
{L_{\bf E}}(V)-\frac{n-2}{2}{\Gamma_{\bf E}}(V)=\frac{\alpha^2n}{2\varphi^2}, 
$$
and the CKN inequality~\eqref{1} takes the form~\eqref{179b}, as announced.  
\medskip

Next, we prove that Sobolev's inequality in the CKN Euclidean space implies the Sobolev inequality on the CKN hyperbolic space. We mimic the previous proof. Define the function $\psi$ on the punctured open unit ball $\mathbf B\setminus\{0\}\subset\R^d$ by 
$$
\psi(x)=\frac{1-|x|^{2\alpha}}2,\quad x\in\mathbf B\setminus\{0\}.
$$
Then, on $\mathbf B\setminus\{0\}$, 
$$
{\Gamma_{\bf H}}=\psi^2{\Gamma_{\bf E}}\quad\text{and}\quad{\mu_{\bf H}}=\psi^{-n}{\mu_{\bf E}}. 
$$
Apply the CKN inequality~\eqref{1} to the function $f=\psi^{\frac{2-n}{2}}h$, where $h\in \mathcal C^\infty_c(\B\setminus\{0\})$.  Again,  we get 
$$
\int f^pd{\mu_{\bf E}}=\int f^{\frac{2n}{n-2}}d{\mu_{\bf E}}=\int h^pd{\mu_{\bf H}}. 
$$
and~\eqref{1} becomes, with $U=\log\psi$, 
$$
\PAR{\int h^pd{\mu_{\bf H}}}^{2/p}\leq C\left(\int {\Gamma_{\bf H}}(h) d{\mu_{\bf H}}+\frac{n-2}{2}
\int \PAR{{ L_{\bf E}}(U)-\frac{n-2}{2}{\Gamma_{\bf E}}(U)}\psi^2h^2 d{\mu_{\bf H}}\right).
$$
We obtain
$$
{L_{\bf E}}(U)-\frac{n-2}{2}{\Gamma_{\bf E}}(U)=-\frac{\alpha^2n}{2\psi^2}, 
$$
and so~\eqref{179c}, as claimed.
\proofend\\

\section{Sobolev's inequality for the spherical CKN model}
\label{sec-3}

This section is devoted to the proof of the optimal Sobolev inequality for the spherical CKN space (Theorem~\ref{thm-1}) under the Felli-Schneider condition~\eqref{fs}. 
It is convenient to introduce spherical coordinates $\R^d\setminus\{0\}\ni x=r\theta$ with $r>0$ and $\theta\in\Sp^{d-1}$. The Sobolev inequality on the CKN sphere~\eqref{355} then takes the form 
\begin{multline*}
\PAR{\int_{(0,\infty)\times\Sp^{d-1}} \vert v\vert^p \frac{r^{d-1-pb}}{(1+r^{2\alpha})^{n}}drdV_{\bf S^{d-1}}}^{2/p}\leq C\int_{(0,\infty)\times\Sp^{d-1}} {\SBRA{(\partial_r v)^2+\frac{1}{r^2}|\nabla_\theta v|^2 }}{\frac{r^{2(\alpha-1)+d-1-pb}}{(1+r^{2\alpha})^{n-2}}}drdV_{\bf S^{d-1}}\\
+4Z^{-\frac2n}\int_{(0,\infty)\times\Sp^{d-1}}v^2 \frac{r^{d-1-pb}}{(1+r^{2\alpha})^{n}}drdV_{\bf S^{d-1}},
\end{multline*}
where $|\nabla_\theta v|$ is the Riemannian length of the Riemannian gradient $\nabla_\theta v$ on $\Sp^{d-1}$ and $dV_{\bf S^{d-1}}$ is the associated Riemannian volume.  Using the  change of variable $(0,\infty)\ni r=e^s$, with $s\in\R$, the inequality becomes (with a different constant $C$), 
\begin{multline}
\label{456}
\PAR{\int_{\R\times\Sp^{d-1}} \vert v\vert^p\cosh(\alpha s)^{-n}dsdV_{\bf S^{d-1}}}^{2/p}\!\!\!\!\leq\\ 
C\int_{\R\times\Sp^{d-1}} \!\!\!\!\!\!\!\!\!\!{\SBRA{(\partial_s v)^2+|\nabla_\theta v|^2}}\cosh(\alpha s)^{2-n}dsdV_{\bf S^{d-1}}
	+Z^{-\frac2n}\!\!\int_{\R\times\Sp^{d-1}} \!\!\!\!\!\!\!\!\!\!\!\!v^2\cosh(\alpha s)^{-n}dsdV_{\bf S^{d-1}},
\end{multline}
where we used the fact that  $d-n\alpha -pb=0$, see~\eqref{99}. 
 This new chart is often called the Emden-Fowler transformation, as suggested in~\cite{florin-wang2001,DolEstLos16}.
In other words, in the cylindrical chart $(s,\theta)\in\R\times\Sp^{d-1}$, the spherical CKN space takes a new and nice form. Notice that the space remains the same, it is only written in a new chart. More precisely, letting 
\begin{equation}
  \label{3540}
  \phi(s)=\cosh(\alpha s), \quad s\in\R,
\end{equation}
 the metric becomes (with the upper indices)
\begin{equation}
  \label{756}
  {\mathfrak g_{\bf S}}=\phi^2\mathfrak h=e^{2\tau_{\bf S}}\mathfrak h,
\end{equation}
where $\tau_{\bf S}=\log\phi$ and 	$\mathfrak h$ is the standard product metric\footnote{on the cotangent space of $(0,\infty)\times \Sp^{d-1}$, where $\Sp^{d-1}$ is viewed in a given chart } on $(0,\infty)\times \Sp^{d-1}$, represented by the $d$-dimensional matrix
$$
\left(
\begin{array}{cc}
  1 & 0\\
  0 & G_\theta
\end{array}
\right),
$$
where \(G_\theta\) is the matrix of $\mathfrak g_\theta$ in the chart \((s,\theta)\), and \(\mathfrak g_\theta\) is the round metric of $\Sp^{d-1}$. For convenience, in Lemma~\ref{lem-2} and its proof, as well as the proof of Proposition~\ref{prop-100}, we will abuse the notations and identify the tensors with their coordinates in the chart \((s,\theta)\), since it will be the only chart used in all the calculations.
The carré du champ operator takes the form   
$$
\Gamma_{\bf S}(f)= \phi^2\SBRA{(\partial_s f)^2+|\nabla_\theta f|^2} =\phi^2\SBRA{(\partial_s f)^2+\Gamma^\theta(f)},
$$
where $\Gamma^\theta(f)=|\nabla_\theta f|^2$ is the carré du champ operator associated to the the Laplace-Beltrami operator $\Delta_{\theta}$ on $\Sp^{d-1}$.  
The Riemannian volume becomes
$
d{V_{{\mathfrak g_{\bf S}}}}=\phi^{-d}{dsdV_{\bf S^{d-1}}}
$
and the reference measure (not normalized measure)
$$
d{\mu_{\bf S}}=\phi^{-n}{dsdV_{\bf S^{d-1}}}.
$$
The corresponding weight $W_{\bf S}$ is defined by $d{\mu_{\bf S}}=e^{-W_{\bf S}}d{V_{{\mathfrak g_{\bf S}}}}$, so that
$$
{W_{\bf S}}=(n-d)\log \phi.
$$
Finally, the associated generator takes the  pleasant form 
$$
{L_{\bf S}}(f)=\phi^2\Big[\partial_{ss}f+(2-n)\frac{\phi'}{\phi}\partial_sf+\Delta_\theta f\Big].
$$
Taking advantage of this chart, let us begin by proving that the spherical CKN space satisfies the $CD(\rho,n)$ condition whenever condition~\eqref{bdgz} holds:

\begin{eproofof}{Proposition~\ref{prop-100}}
From~\cite[Sec.~C6]{bgl-book}, the generator ${L_{\bf S}}$ satisfies a $CD(\rho, n)$ condition (with $n>d$)  if and only if, as a symmetric tensor (with lower indices),
$$
Ric({L_{\bf S}})-\rho{ g_{\bf S}}\geq \frac{1}{n-d}\nabla^{{\mathfrak g_{\bf S}}}{W_{\bf S}}\otimes \nabla^{{\mathfrak g_{\bf S}}}{W_{\bf S}}.
$$ 
 Let us remark that, since \({\mathfrak g_{\bf S}}= \phi^2\mathfrak h\) (with upper indices), the corresponding metric tensors (with lower indices) satisfy 
$$
{ g_{\bf S}}=\frac{{h}}{\phi^2}.
$$

Compute first the r.h.s. of the above inequality. From the definition of ${W_{\bf S}}$, we have
$$
\frac{\nabla^{{\mathfrak g_{\bf S}}}{W_{\bf S}}\otimes \nabla^{{\mathfrak g_{\bf S}}}{W_{\bf S}}}{n-d}=(n-d)\PAR{\frac{\phi'}{\phi}}^2 J,
$$
where  $J$ is the $d$-dimensional matrix with all entries equal to zero but the first i.e. $J_{ij}=\delta_{i1}\delta_{j1}$ or more visually, letting $H$ (resp. $G_\theta$) be the matrix representing the standard product metric $h$ (resp. $g_\theta$) in the coordinates $(s,\theta)$ (resp. $\theta$), 
\begin{equation}
  \label{758}
  J=
  \begin{pmatrix}
    1 & 0\\
    0 & 0
  \end{pmatrix}
  =
  H-
  \begin{pmatrix}
    0 & 0\\
    0 & G_\theta
  \end{pmatrix}.
\end{equation}
Now, applying formula~\eqref{230} in Lemma~\ref{lem-2} below,  we get
$$
Ric({L_{\bf S}})-\rho{ g_{\bf S}}- \frac{1}{n-d}\nabla^{{\mathfrak g_{\bf S}}}{W_{\bf S}}\otimes \nabla^{{\mathfrak g_{\bf S}}}{W_{\bf S}}=- (n-d)\alpha^2 H+(d-2)(1-\alpha^2)(H- J)+ (n-d)\alpha^2J,
$$
where, again, we conflate tensors and their matrices in the chart \((s,\theta)\). Hence,  
$$
Ric({L_{\bf S}})-\rho{ g_{\bf S}}- \frac{1}{n-d}\nabla^{{\mathfrak g_{\bf S}}}{W_{\bf S}}\otimes \nabla^{{\mathfrak g_{\bf S}}}{W_{\bf S}}=\mathfrak B_{DGZ}
\begin{pmatrix}
    0 & 0\\
    0 & G_\theta
\end{pmatrix}
$$
with the constant  $\mathfrak B_{DGZ}$ is defined given in~\eqref{951},
and so,  ${L_{\bf S}}$ satisfies the curvature-dimension condition $CD(\rho,n)$ if and only if  $\mathfrak B_{DGZ}\geq0$. 
\end{eproofof}
\begin{erem}Since the matrix $H- J$ depends only on the variable $\theta$, when we restrict to functions depending on the variable $s$ only, the corresponding model {\it always} satisfies the $CD(\rho,n)$ condition, regardless of the sign of $\mathfrak B_{DGZ}$. 
\end{erem}
In the above proof, we made strong use of the following lemma.
\begin{elem}[Computation of $Ric({L_{\bf S}})$]
\label{lem-2}
We have the following formulae 
\begin{equation}
  \label{21}
  Ric_{{ g_{\bf S}}}=\frac{(d-1)\alpha^2}{\phi^2}
  \begin{pmatrix}
    1 & 0 \\
    0 & 0
  \end{pmatrix}
  +\frac{1}{\phi^2}\SBRA{(d-2)(1-\alpha^2)\phi^2+(d-1)\alpha^2}
    \begin{pmatrix}
    0 & 0\\
    0 & G_\theta
  \end{pmatrix},
\end{equation}
and
\begin{equation}
  \label{23}
  \nabla\nabla^{{\mathfrak  g_{\bf S}}}\, W_{\bf S}=(n-d)\alpha^2
  \begin{pmatrix}
    1 & 0 \\
    0 & 0
  \end{pmatrix}
  + (n-d)\alpha^2\frac{1-\phi^2}{\phi^2}
  \begin{pmatrix}
    0 & 0\\
    0 & G_\theta
  \end{pmatrix},
\end{equation}
where \(G_\theta\) is the matrix of round metric on the sphere $\Sp^{d-1}$. With the constant  $\mathfrak B_{DGZ}$ (given in~\eqref{951}), we obtain
\begin{equation}
  \label{230}
  Ric({L_{\bf S}})=Ric_{{ g_{\bf S}}}+\nabla\nabla^{{\mathfrak  g_{\bf S}} }{W_{\bf S}}=\frac{\alpha^2}{\phi^2}\SBRA{d-1+\phi^2(n-d)}
  \begin{pmatrix}
    1 & 0 \\
    0 & 0
  \end{pmatrix}
  +\frac{1}{\phi^2}\SBRA{\alpha^2(n-1)+\phi^2\mathfrak B_{DGZ}}
  \begin{pmatrix}
    0 & 0\\
    0 & G_\theta
  \end{pmatrix}.
\end{equation}
\end{elem}
\begin{eproof}
Let us start with $Ric_{{ g_{\bf S}}}$, which is simply the Ricci tensor of the metric ${\mathfrak g_{\bf S}}$.  Since 
$
{\mathfrak g_{\bf S}}=e^{2\tau_{\bf S}}\mathfrak h
$
is conformal to $\mathfrak h$, we may apply~\eqref{120} in the appendix to get (with lower indices)
\begin{equation}
\label{20}
Ric_{{ g_{\bf S}}}=Ric_{h}+(\Delta_{\mathfrak h}\tau_{\bf S}){h}+(d-2)(\nabla\nabla^{\mathfrak h}\tau_{\bf S}+\nabla^{\mathfrak h}\tau_{\bf S} \odot_{\mathfrak h} \nabla^{\mathfrak h}\tau_{\bf S}-\nabla^{\mathfrak h}\tau_{\bf S}\cdot_{\mathfrak h}\nabla^{\mathfrak h}\tau_{\bf S}\, { h}).
\end{equation}
Since $Ric_{g_\theta}=(d-2)g_\theta$, we have 
$$
Ric_{h}=(d-2)(H-J)=(d-2)  \begin{pmatrix}
    0 & 0\\
    0 & G_\theta
  \end{pmatrix}.
$$
Since $\phi$ depends only on the variable $s$, we have 
\begin{align*}
  &\Delta_{\mathfrak h}(\tau_{\bf S})=\tau_{\bf S}''=\frac{\phi''}{\phi}-\PAR{\frac{\phi'}{\phi}}^2=\alpha^2-\PAR{\frac{\phi'}{\phi}}^2, \\
  &\nabla^{\mathfrak h}\tau_{\bf S} \odot_{\mathfrak h} \nabla^{\mathfrak h}\tau_{\bf S} = \PAR{\frac{\phi'}{\phi}}^2J,\\
  &\nabla^{\mathfrak h}\tau_{\bf S}=\frac{\phi'}{\phi}
  \begin{pmatrix}
    1 \\
    0
  \end{pmatrix},
\end{align*}
and 
 $$
 (\nabla^{\mathfrak h}\tau_{\bf S}\cdot_{\mathfrak h}\nabla^{\mathfrak h}\tau_{\bf S})\, { h}=\PAR{\frac{\phi'}{\phi}}^2 H.
 $$
Collecting the four terms and using~\eqref{20}, we get 
$$
Ric_{g}= H\SBRA{d-2+\alpha^2-(d-1)\PAR{\frac{\phi'}{\phi}}^2}+J(d-2)(\alpha^2-1).
$$
Since $\phi'^2=\alpha^2(\phi^2-1)$, the equation  can be written, 
$$
Ric_{ g}= H\frac{(d-1)\alpha^2}{\phi^2}+(H-J)(d-2)(1-\alpha^2),
$$
which is the desired result. 

Let us now compute $\nabla\nabla^{{\mathfrak  g_{\bf S}}} W_{\bf S}$, the Hessian with respect to the metric ${\mathfrak g_{\bf S}} $. We have (see~\eqref{130}), 
$$
\nabla\nabla^{{\mathfrak  g_{\bf S}}}\, W_{\bf S}={\nabla\nabla}^{{\mathfrak  h}} W_{\bf S}+2\nabla W_{\bf S} \odot_{{\mathfrak  h}}\nabla^{{\mathfrak  h}}\tau_{\bf S}-(\nabla^{{\mathfrak  h}} W_{\bf S}\cdot_{{\mathfrak  h}}\nabla^{{\mathfrak  h}} \tau_{\bf S})\,{{ h}}.
$$
Since ${W_{\bf S}}$ depends only on the variable $s$, we easily get that
$$
\nabla\nabla^{{\mathfrak  g_{\bf S}}}\, W_{\bf S}=J(n-d)\SBRA{\alpha^2-\PAR{\frac{\phi'}{\phi}}^2}+2J (n-d)\PAR{\frac{\phi'}{\phi}}^2 - H(n-d)\PAR{\frac{\phi'}{\phi}}^2,
$$
which is the expected result. 
\end{eproof}

\begin{erem} As an immediate consequence of Proposition~\ref{prop-100}, the fact that $L_{\bf S}$ is essentially self-adjoint  when $n\ge 3$ (see~\cite[Thm.~3.12]{Ketterer}) and Theorem~\ref{th:B},  
we see that Sobolev's inequality~\eqref{355} holds (and so Poincaré's inequality~\eqref{890} too), as soon as~\eqref{bdgz} holds. Also note that $f=\frac{\phi'}{\phi}$, seen as a function of the first of the cylindrical coordinates $(s,\theta)$, solves
$$
-{L_{\bf S}}f=n\alpha^2 f
$$
and so equality in Poincaré's inequality~\eqref{890} is achieved by $f$. In particular, the constant in Sobolev's inequality~\eqref{355} is optimal.
\end{erem}
In fact, one can do better and prove optimal inequalities in the optimal range of parameters given by the Felli-Schneider condition, as we describe next. 
The first crucial step consists in proving the  following {\it weaker} integrated forms of the curvature-dimension condition~\eqref{cdrhoi}.
\begin{eprop}
\label{crucial} Let $(a,b)\in\Theta_{FS}$. In cylindrical coordinates $(s,\theta)$, for any $s\in(0,\infty)$ and any smooth positive function $f$ on $(0,\infty)\times \Sp^{d-1}$, there holds 
\begin{equation}
\int \PAR{{\Gamma^{\bf S}_2}(f)-\rho{\Gamma_{\bf S}}(f)-\frac{1}{n}({L_{\bf S}} f)^2}f^{1-n}d{V_{\bf S^{d-1}}}\geq 0 
    \label{icd}
\end{equation}
and
\begin{equation}\label{icd2}
\int \PAR{{\Gamma^{\bf S}_2}(f)-\rho{\Gamma_{\bf S}}(f)-\frac{1}{n}({L_{\bf S}} f)^2}d{V_{\bf S^{d-1}}}\geq 0, 
\end{equation}
where $d{V_{\bf S^{d-1}}}$ is the standard volume on the sphere $\bf S^{d-1}$.
\end{eprop}
We establish Proposition~\ref{crucial} through a series of lemmas. First,
\begin{elem}[${\Gamma^{\bf S}_2}$ in the cylindrical chart]
\label{lem-147}
	Let $(a,b)\in\Theta$. In cylindrical coordinates, we have for any smooth function $f$ on $(0,\infty)\times \Sp^{d-1}$ 
	\begin{multline}
	\label{564}
	\frac{{\Gamma^{\bf S}_2}(f)}{\phi^4}=
	(\partial_{ss}f)^2+||\nabla\nabla_\theta f||^2+2\Gamma^\theta(\partial_s f)+2\frac{\phi'}{\phi} \partial_{ss}f\partial_s f+4\frac{\phi'}{\phi}\Gamma^\theta(\partial_sf,f)	-2\frac{\phi'}{\phi}\partial_s f\Delta_\theta f\\
	(\partial_sf)^2\SBRA{d\PAR{\frac{\phi'}{\phi}}^2+\alpha^2\PAR{\frac{d-1}{\phi^2}+n-d}}+\Gamma^\theta(f)\PAR{2\PAR{\frac{\phi'}{\phi}}^2+\alpha^2\frac{n-1}{\phi^2}+\mathfrak B_{DGZ}},
	\end{multline}
	where $||\nabla\nabla_\theta f||^2$ is the Hilbert-Schmidt norm with respect to the variable $\theta$, $\Gamma^\theta(f)=|\nabla_\theta f|^2$ the carré du champ operator associated to $\Delta_{\theta}$ and the function  $\phi$ has been defined in~\eqref{3540}.
\end{elem}
\begin{eproof}
We can use the definition of the $\Gamma_2$ operator to prove~\eqref{564}. But, since the Ricci curvature of ${L_{\bf S}}$ has been computed in Lemma~\ref{lem-2}, we prefer to use the following Bochner-Lichnerowicz formula,
	\begin{equation}
	\label{235}
	{\Gamma^{\bf S}_2}(f)=Ric({L_{\bf S}})(\nabla f,\nabla f)+||\nabla\nabla^{{\mathfrak g_{\bf S}}} f||^2,
	\end{equation}
	where $||\nabla\nabla^{{\mathfrak g_{\bf S}}}f||^2$ is the Hilbert-Schmidt norm of the Hessian of $f$ with respect to the metric ${\mathfrak g_{\bf S}}$ (see for instance~\cite[P.~71]{bgl-book}).  From Lemma~\ref{lem-2}, equation~\eqref{230}, we have first 
	$$
	\frac{Ric({L_{\bf S}})(\nabla f,\nabla f)}{\phi^4}=\frac{(\partial_sf)^2}{\phi^2}\alpha^2\SBRA{d-1+\phi^2(n-d)}+\frac{\Gamma^\theta(f)}{\phi^2}\SBRA{\alpha^2(n-1)+\phi^2\mathfrak B_{DGZ}}.
	$$

	It remains to compute $||\nabla\nabla^{{\mathfrak g_{\bf S}}}f||^2$. 
	From~\eqref{756} we have  $	{\mathfrak g_{\bf S}}=\phi^2\mathfrak h=e^{2\tau_{\bf S}}\mathfrak h$ and so we may apply formula~\eqref{131} to get
$$
	\frac{||{\nabla\nabla}^{{\mathfrak g_{\bf S}}}  f||^2}{\phi^4}=||\nabla\nabla^{\mathfrak h} f||^2+2\Gamma^{\mathfrak h}(\tau_{\bf S},\Gamma^{\mathfrak  g}(f))+2\Gamma^{\mathfrak  h}(f)\Gamma^{\mathfrak  h}(\tau_{\bf S})+
	(d-2)\Gamma^{\mathfrak  h}(f,\tau_{\bf S})^2-2\Delta_{\mathfrak h} f\Gamma^{\mathfrak  h}(f,\tau_{\bf S}).
	$$
	Since $\mathfrak h$ is the standard metric product and $\tau$ depends only on the variable $s$, we have 
	$$
	||\nabla\nabla^{\mathfrak h} f||^2=(\partial_{ss} f)^2+||\nabla\nabla_\theta f||^2+2\Gamma^\theta(\partial_s f),
	$$ 
	$$
\Gamma^{\mathfrak h}(\tau_{\bf S},\Gamma^{\mathfrak  g}(f))=2\frac{\phi'}{\phi}\partial_{ss}f\partial_s f+	2\frac{\phi'}{\phi}\Gamma^\theta(\partial_s f,f),
	$$
	$\Gamma^{\mathfrak  h}(f)=(\partial_sf)^2+\Gamma^\theta(f)$, 	$\Gamma^{\mathfrak  h}(\tau_{\bf S})=\PAR{\frac{\phi'}{\phi}}^2$, $\Gamma^{\mathfrak  h}(f,\tau_{\bf S})=\frac{\phi'}{\phi}\partial_sf$ and 
	$\Delta_{\mathfrak h} f=\partial_{ss} f+\Delta_{\theta}f$.
	
	Collecting all the terms, we get 
	 \begin{multline*}
	 \frac{||{\nabla\nabla}^{{\mathfrak g_{\bf S}}}  f||^2}{\phi^4}=(\partial_{ss} f)^2+||\nabla\nabla_\theta f||^2+2\Gamma^\theta(\partial_s f)+4\frac{\phi'}{\phi}\partial_{ss}f\partial_s f+4\frac{\phi'}{\phi}\Gamma^\theta(\partial_s f,f)\\
	 +2\PAR{\frac{\phi'}{\phi}}^2\SBRA{(\partial_sf)^2+\Gamma^\theta(f)}+
	 (d-2)\PAR{\frac{\phi'}{\phi}}^2(\partial_sf)^2-2\PAR{\partial_{ss} f+\Delta_{\theta}f}\frac{\phi'}{\phi}\partial_sf,
	 \end{multline*}
	 that is
	 \begin{multline*}
	 \frac{||{\nabla\nabla}^{{\mathfrak g_{\bf S}}}  f||^2}{\phi^4}=(\partial_{ss} f)^2+||\nabla\nabla_\theta f||^2+2\Gamma^\theta(\partial_s f)+2\frac{\phi'}{\phi}\partial_{ss}f\partial_s f+4\frac{\phi'}{\phi}\Gamma^\theta(\partial_s f,f)\\
	 +2\PAR{\frac{\phi'}{\phi}}^2\Gamma^\theta(f)+
	 d\PAR{\frac{\phi'}{\phi}}^2(\partial_sf)^2-2\Delta_{\theta}f\frac{\phi'}{\phi}\partial_sf.
	 \end{multline*}
	 Finally, by using~\eqref{235}, we get the expected formula~\eqref{564}. 
\end{eproof}
We restate the above lemma in the following more compact formulation.
\begin{elem}
\label{lem-12}
       In the cylindrical chart, for any smooth function $f$ on $(0,\infty)\times \Sp^{d-1}$,  
	\begin{multline}
	\label{554}
	\frac{1}{\phi^4}\PAR{{\Gamma^{\bf S}_2}(f)-\rho{\Gamma_{\bf S}}(f)-\frac{1}{n}({L_{\bf S}} f)^2}=
	\frac{n-1}{n}\PAR{\partial_{ss}f+2\frac{\phi'}{\phi}\partial_s f -\frac{1}{n-1}\Delta_\theta f}^2+||\nabla\nabla^\theta f||^2-\frac{1}{n-1}(\Delta_\theta f)^2\\
	+2\Gamma^\theta\big(\partial_s f+\frac{\phi'}{\phi}f\big)+\Gamma^\theta(f)\mathfrak B_{DGZ}.
	\end{multline}
\end{elem}

\begin{eproof}
In the cylindrical chart, the generator takes the following form, for a smooth function $f$:
$$
{L_{\bf S}}(f)=\phi^2\Big[\partial_{ss}f+(2-n)\frac{\phi'}{\phi}\partial_sf+\Delta_\theta f\Big],
$$
and  from Lemma~\ref{lem-147} (formula~\eqref{564}), we obtain 
 \begin{multline*}
\frac{1}{\phi^4}\SBRA{{\Gamma^{\bf S}_2}(f)-\rho{\Gamma_{\bf S}}(f)-\frac{1}{n}({L_{\bf S}} f)^2}=
\frac{n-1}{n}(\partial_{ss}f)^2+||\nabla\nabla_\theta f||^2-\frac{1}{n}(\Delta_\theta f)^2 +4\frac{\phi'}{\phi} \frac{n-1}{n}\partial_{ss}f\partial_s f-\frac{4}{n}\frac{\phi'}{\phi}\partial_s f\Delta_\theta f\\
-\frac{2}{n}\partial_{ss} f\Delta_\theta f
+2\Gamma^\theta(\partial_sf)+4\frac{\phi'}{\phi}\Gamma^\theta(\partial_sf,f)\\
+(\partial_sf)^2\SBRA{\PAR{\frac{\phi'}{\phi}}^2\PAR{d-\frac{(n-2)^2}{n}}+(n-d)\alpha^2\frac{\phi^2-1}{\phi^2}}+\Gamma^\theta(f)\PAR{2\PAR{\frac{\phi'}{\phi}}^2+\mathfrak B_{DGZ}}.
\end{multline*}
Since $\phi'^2=\alpha^2(\phi^2-1)$, we get 
 \begin{multline*}
	\frac{1}{\phi^4}\SBRA{{\Gamma^{\bf S}_2}(f)-\rho{\Gamma_{\bf S}}(f)-\frac{1}{n}({L_{\bf S}} f)^2}=
	\frac{n-1}{n}(\partial_{ss}f)^2+||\nabla\nabla_\theta f||^2-\frac{1}{n}(\Delta_\theta f)^2 +4\frac{\phi'}{\phi} \frac{n-1}{n}\partial_{ss}f\partial_s f-\frac{4}{n}\frac{\phi'}{\phi}\partial_s f\Delta_\theta f\\
	+4\frac{n-1}{n}\PAR{\frac{\phi'}{\phi}}^2(\partial_sf)^2-\frac{2}{n}\partial_{ss} f\Delta_\theta f
	+\Gamma^\theta(f)\PAR{2\PAR{\frac{\phi'}{\phi}}^2+\mathfrak B_{DGZ}}+2\Gamma^\theta(\partial_sf)+4\frac{\phi'}{\phi}\Gamma^\theta(\partial_sf,f).
	\end{multline*}
Formula~\eqref{554} follows then  easily since $\phi$ depends only on the variable $s$. 	
\end{eproof}

\begin{erem}
By the Cauchy-Schwarz inequality, 
$$
||\nabla\nabla_\theta f||^2\geq \frac{1}{d-1}(\Delta_\theta f)^2,
$$
and so
$$
||\nabla\nabla_\theta f||^2-\frac{1}{n-1}(\Delta_\theta f)^2\geq\frac{n-d}{(d-1)(n-1)} (\Delta_\theta f)^2\geq0
$$
since $n\geq d$. We recover from \eqref{554} that, under the condition $\mathfrak B_{DGZ}\geq 0$, the generator ${L_{\bf S}}$ satisfies the $CD(\rho,n)$ curvature-dimension condition. 
\end{erem}
The next ingredient is the following inequality valid on the sphere ${\bf S}^{d-1}$ (or any smooth weighted manifold satisfying the $CD(d-2,d-1)$ condition).
\begin{elem}
\label{lem-23}
For any smooth positive function $f$  on $\Sp^{d-1}$,
\begin{equation}
\label{687}
\int \Gamma^\theta_2(f)f^{1-n}d{V_{\bf S^{d-1}}}\geq (d-1)\int \Gamma^\theta(f)f^{1-n}d{V_{\bf S^{d-1}}}+A\int \frac{\Gamma^\theta(f)^2}{f^2}f^{1-n}d{V_{\bf S^{d-1}}},
\end{equation}
where 
\begin{equation}
\label{187}
A=\frac{n-1}{4(d+1)^2}(n(4d-5)+3(4d+7)).
\end{equation}
In particular, 
\begin{equation}
\label{215}
\int \Gamma^\theta_2(f)f^{1-n}d{V_{\bf S^{d-1}}}\geq (d-1)\int \Gamma^\theta(f)f^{1-n}d{V_{\bf S^{d-1}}}.
\end{equation}

\end{elem}
\begin{eproof}
The operator $\Delta_\theta$ is the Laplace-Beltrami operator on the $(d-1)$-dimensional sphere, therefore, it satisfies the $CD(d-2,d-1)$ condition. Moreover,   $Ric_{\mathfrak g_\theta}$ (the Ricci tensor of $\Delta_\theta$) satisfies  
\begin{equation}
\label{159}
Ric_{\mathfrak g_\theta}(\nabla_\theta f,\nabla_\theta f)=(d-2)\Gamma^\theta(f).
\end{equation} 
In~\cite[p.~767]{gentil-zug}, it is proved that under the $CD(K,m)$ condition, for a general operator $L$ associated to the measure  $\mu$, and the operators $\Gamma$ and $\Gamma_2$, one has for any real parameters $q,\chi$,
  $$
  \int h^q\Gamma_2(h)d\mu\geq \frac{K m}{m-1}\int h^q\Gamma(h)d\mu+ \int \SBRA{Ah^{q-2}\Gamma( h)^2+B h^{q-1} \Gamma(h,\Gamma(h))}d\mu
  $$
  where
  $$
  \left\{
  \begin{aligned}
    A&=\frac{q(q-1)}{m-1}-\chi^2-2\chi\frac{q-1}{m-1},\\
    B&=\frac{1}{m-1}\left(\frac{3q}{2}-\chi(m+2)\right). 
  \end{aligned}\right.
  $$
Apply the previous inequality to our operator  $\Delta_\theta$ with parameters  $q=1-n$, $K=d-2$, $m=d-1$ and $\chi=\frac{3q}{2(m+2)}$ so that $B=0$.  We obtain 
\begin{equation*}
  \int \Gamma^\theta_2(f)f^{1-n}dV_{\bf S^{d-1}}
   \geq (d-1)\int_{\Sp^{d-1}} \Gamma^\theta(f)f^{1-n} dV_{\bf S^{d-1}} +A\int_{\Sp^{d-1}}  \frac{\Gamma^\theta(f)^2}{f^2}f^{1-n}dV_{\bf S^{d-1}} 
\end{equation*}
where $A$ is given by~\eqref{187} after a straightforward computation. In particular, $A\geq0$ and the estimate~\eqref{215} follows. 
\end{eproof}
We can now turn to the
\begin{eproofof}{Proposition~\ref{crucial}}
From Lemma~\ref{lem-12}, we have 
$$
\int \PAR{{\Gamma^{\bf S}_2}(f)-\rho{\Gamma_{\bf S}}(f)-\frac{1}{n}({L_{\bf S}} f)^2}f^{1-n}d{V_{\bf S^{d-1}}}\geq \phi^4\int \PAR{||\nabla\nabla^\theta f||^2-\frac{1}{n-1}(\Delta_\theta f)^2+\Gamma^\theta(f)\mathfrak B_{DGZ}}f^{1-n}d{V_{\bf S^{d-1}}}.
$$
Now,
$$
||\nabla\nabla_\theta f||^2-\frac{1}{n-1}(\Delta_\theta f)^2= \frac{n-d}{n-1}||\nabla\nabla_\theta f||^2+\frac{d-1}{n-1}||\nabla\nabla_\theta f||^2-\frac{1}{n-1}(\Delta_\theta f)^2\geq \frac{n-d}{n-1}||\nabla\nabla_\theta f||^2,
$$ 
where we used the Cauchy-Schwarz inequality to infer that
\begin{equation}\label{opt1}
||\nabla\nabla_\theta f||^2\geq \frac{1}{d-1}(\Delta_{\theta} f)^2.
\end{equation}
So,
$$
\int \PAR{{\Gamma^{\bf S}_2}(f)-\rho{\Gamma_{\bf S}}(f)-\frac{1}{n}({L_{\bf S}} f)^2}f^{1-n}d{V_{\bf S^{d-1}}}\geq \phi^4\int \PAR{\frac{n-d}{n-1}||\nabla\nabla^\theta f||^2+\Gamma^\theta(f)\mathfrak B_{DGZ}}f^{1-n}d{V_{\bf S^{d-1}}}.
$$

Since, from~\eqref{235} and~\eqref{159},  
$$
||\nabla\nabla^\theta f||^2=\Gamma^\theta_2(f)-Ric_{\Sp^{d-1}}(\nabla_\theta f,\nabla_\theta f)=\Gamma^\theta_2(f)-(d-2)\Gamma^\theta(f),
$$
the inequality becomes 
\begin{multline*}
\phi^{-4}\int \PAR{{\Gamma^{\bf S}_2}(f)-\rho{\Gamma_{\bf S}}(f)-\frac{1}{n}({L_{\bf S}} f)^2}f^{1-n}d{V_{\bf S^{d-1}}}\geq \\
\frac{n-d}{n-1}\int \Gamma^\theta_2(f)f^{1-n}d{V_{\bf S}^{d-1}}+\SBRA{\frac{n-d}{n-1}(d-2)+\mathfrak B_{DGZ}}\int \Gamma^\theta(f)f^{1-n}d{V_{\bf S^{d-1}}}.
\end{multline*}

Using the estimate~\eqref{215} in Lemma~\ref{lem-23},  we get 
\begin{equation}\label{optr}
\phi^{-4}\int \PAR{{\Gamma^{\bf S}_2}(f)-\rho{\Gamma_{\bf S}}(f)-\frac{1}{n}({L_{\bf S}} f)^2}f^{1-n}d{V_{\bf S^{d-1}}}\geq \SBRA{\frac{n-d}{n-1}+\mathfrak B_{DGZ}}\int\Gamma^\theta(f)f^{1-n}d{V_{\bf S^{d-1}}}.
\end{equation}
That is, 
\begin{equation}\label{opt2}
\frac{n-d}{n-1}+\mathfrak B_{DGZ}\geq 0
\end{equation}
implies \eqref{icd}. The proof of inequality \eqref{icd2} is almost identical, except that instead of \eqref{687}, one uses 
$$
\int \Gamma^\theta_2(f)d{V_{\bf S^{d-1}}}\geq (d-1)\int \Gamma^\theta(f)d{V_{\bf S^{d-1}}},
$$
which itself holds thanks to the Cauchy-Schwarz inequality \eqref{opt1}, Bochner's formula \eqref{235} and the identity $\int \Gamma^\theta_2(f)d{V_{\bf S^{d-1}}}=\int (\Delta_\theta f)^2d{V_{\bf S^{d-1}}}$.
\end{eproofof}
Now that the integrated curvature-dimension is established, we can turn to the proof of Sobolev's inequality.
\begin{eproofof}{Theorem~\ref{thm-1}}
Fix $q\in[1,p)$. By the Caffarelli-Kohn-Nirenberg inequality \eqref{1} and Theorem \ref{confinvckn}, Sobolev's inequality holds on the CKN spherical space in the form \eqref{179b}. By Proposition 6.2.2 in \cite{bgl-book}, \eqref{179b} and Poincaré's inequality \eqref{cor:poincare} imply the following tight form of Sobolev's inequality: 
\begin{equation}\label{eq:sobno}
\PAR{\int \ABS{v}^q d\mu }^{2/q} \le A\int\Gamma_{\bf S}(v) d\mu + \int v^2d\mu,
\end{equation}
for some $A\in\R_+^*$ where $\mu=\frac1Z{\mu_{\bf S}}$ is the normalized measure and $v\in H^1_0(\mu_{\bf S})$. 
Given $A\in\R_+^*$, consider the minimization problem
$$
I(A)=\inf \left\{
A \int\Gamma_{\bf S}(v) d\mu + \int v^2d\mu\;:\; v\in H_0^1(\mu_{\bf S})\;,\; ||v||_{L^q(\mu)}=1
\right\}.
$$
Using $v=1$ as a test function, we see that $I(A)\le 1$. Thus,~\eqref{eq:sobno} holds if and only if $I(A)=1$. Thanks to the Banach-Alaoglu-Bourbaki and Lemma \ref{rellich}, there exists a minimizer $v\in H^1_0(\mu_{\bf S})$ s.t. $||v||_{L^q(\mu)}=1$. By Stampacchia's theorem~\cite{Sta66}, $\vert v\vert$ is also a minimizer, so we may assume that $v\ge 0$ a.e. In addition, a  constant multiple of $v$ (abusively denoted the same below) is a weak solution to
\begin{equation}
\label{10}
-A L_{\bf S} v + v = v^{q-1}\quad\text{in $\R^d\setminus\{0\}$}.
\end{equation}
By standard elliptic regularity (see e.g.~\cite{Heb97}, proof of Theorem 6.2.1, p.~248) $v\in \mathcal C^3(\R^d\setminus\{0\})$ and by the strong maximum principle (see e.g.~\cite{Heb97}, Theorem 5.7.2), $v>0$ in $\R^d\setminus\{0\}$. In addition,  
\begin{equation}\label{lestim}
C\ge v\ge c>0\quad\text{ and }\quad\Gamma_{\bf S}(v)\le C
\end{equation}
for some constants $C,c>0$.
The upper bound on $v$ is obtained by standard Moser iteration (i.e. by multiplying~\eqref{10} by the test function $\min(v,k)^{2\alpha-1}$, where $k\in \N^*$ and $\alpha\ge 1$ and making use of Sobolev's inequality \eqref{non tight sobolev} inductively). For the lower bound on $v$, we apply Proposition 6.3.4 in \cite{bgl-book} and repeat the considerations of p.~312 in the same reference. The upper bound on $\Gamma_{\bf S}(v)$ is more delicate and proved in Lemma~\ref{lem-985}.

Define the pressure function $\Phi=v^{-\frac{q-2}2}$. Then, $\Phi$ solves
\begin{equation}\label{eq:p}
 \Phi {L_{\bf S}} \Phi -\frac {\nu}2{\Gamma_{\bf S}}(\Phi) = -\lambda (\Phi^2-1)\quad\text{in $\R^d\setminus\{0\}$,}
\end{equation}
where $\nu=\frac{2q}{q-2}$ and $\lambda=\frac{q-2}{2A}=\frac2{(\nu-2)A}$. Since $v$ is bounded above and below by positive constants and since $\Gamma_{\bf S}(v)$ is bounded, equation \eqref{eq:p} implies that for every $a\in\R$
\begin{equation}\label{phi2dansl2}
    \Phi^a\in D({L_{\bf S}}) 
\end{equation} 
Multiply  equation~\eqref{eq:p} by \({L_{\bf S}} (\Phi^{1-{\nu}})\) and integrate. Thanks to the integration by parts formula \eqref{IPP},  we find for the right-hand-side
\begin{align*}
  \int {\lambda(\Phi^2-1)}{L_{\bf S}}( \Phi^{1-\nu})
  &=  \lambda\int \Phi^2{L_{\bf S}} (\Phi^{1-\nu})=-\lambda\int {\Gamma_{\bf S}} (\Phi^2, \Phi^{1-\nu} )\\
  &=c\int {\Gamma_{\bf S}}(\Phi)\Phi^{1-\nu}
\end{align*}
where $c=2\lambda(\nu-1)=4\frac{\nu-1}{(\nu-2)A}$ and where integration is understood with respect to the reference measure $\mu=\frac1Z\mu_{\bf S}$.
For the left-hand side, integrations by parts must be dealt with more carefully. By Lemma \ref{kesa} (or since ${L_{\bf S}}$ is essentially self-adjoint by \cite[Thm.~3.12]{Ketterer}), there exists a sequence of radial functions $\zeta_k\in C^\infty_c(\R^d\setminus\{0\})$ such that $\zeta_k\to 1$ in $D({L_{\bf S}})$. By equation \eqref{eq:p}, $\Phi {L_{\bf S}} \Phi -\frac {\nu}2{\Gamma_{\bf S}}(\Phi)$ is bounded. By \eqref{phi2dansl2}, ${L_{\bf S}} (\Phi^{1-{\nu}})\in L^2(\mu_{\bf S})$. So, by dominated convergence, as $k\to+\infty$, 
$$
\int ({L_{\bf S}} (\Phi^{1-{\nu}}))(\Phi {L_{\bf S}}\Phi-\frac{\nu}{2}{\Gamma_{\bf S}}(\Phi)) = \int ({L_{\bf S}} (\Phi^{1-{\nu}}))(\Phi {L_{\bf S}}\Phi-\frac{\nu}{2}{\Gamma_{\bf S}}(\Phi))\zeta_k + o(1)
$$
Since $\zeta_k$ is compactly supported, we may integrate by parts and deduce that
\begin{align}
 \int ({L_{\bf S}} (\Phi^{1-{\nu}}))(\Phi {L_{\bf S}}\Phi-\frac{\nu}{2}{\Gamma_{\bf S}}(\Phi))\zeta_k &= -\int{\Gamma_{\bf S}}((\Phi {L_{\bf S}}\Phi-\frac{\nu}{2}{\Gamma_{\bf S}}(\Phi))\zeta_k,\Phi^{1-\nu}) + o(1)\nonumber\\  
 &= -\int{\Gamma_{\bf S}}((\Phi {L_{\bf S}}\Phi)\zeta_k,\Phi^{1-\nu}) + \frac{\nu}{2}\int{\Gamma_{\bf S}}({\Gamma_{\bf S}}(\Phi))\zeta_k,\Phi^{1-\nu}) + o(1)\nonumber\\
 &=: I+II+o(1)\label{rhs}
\end{align}
Using the product rule to expand derivatives in the first integral, we find
$$
I=-\int \zeta_k {L_{\bf S}}\Phi{\Gamma_{\bf S}}(\Phi,\Phi^{1-\nu})-\int\zeta_k\Phi{\Gamma_{\bf S}}( {L_{\bf S}}\Phi,\Phi^{1-\nu})-\int(\Phi {L_{\bf S}}\Phi){\Gamma_{\bf S}}(\zeta_k,\Phi^{1-\nu})=:I_1+I_2+I_3
$$
Using \eqref{eq:p} and the boundedness of $\Phi$ and ${\Gamma_{\bf S}}(\Phi)$, we find
\begin{equation}\label{I3}
    I_3
\le  C\left(\int {\Gamma_{\bf S}}(\zeta_k)\right)^{1/2}=o(1)
\end{equation}
Next, we deal with $I_1$. Thanks to the product rule for derivatives, we find
\begin{align*}
    I_1&= -\int 
    {\Gamma_{\bf S}}(\Phi,\Phi^{1-\nu}{L_{\bf S}}\Phi \zeta_k)+
    \int\Phi^{1-\nu}\zeta_k {\Gamma_{\bf S}}(\Phi,{L_{\bf S}}\Phi)+
    \int\Phi^{1-\nu}{L_{\bf S}} \Phi {\Gamma_{\bf S}}(\Phi,\zeta_k)\\
    &=-\int 
    {\Gamma_{\bf S}}(\Phi,\Phi^{1-\nu}{L_{\bf S}}\Phi \zeta_k)+
    \int\Phi^{1-\nu}\zeta_k {\Gamma_{\bf S}}(\Phi,{L_{\bf S}}\Phi)+o(1)
\end{align*}
Since $\zeta_k$ has compact support, we may integrate by parts to find that
\begin{equation}\label{I1}
    I_1=-\int\Phi^{1-\nu}({L_{\bf S}}\Phi)^2\zeta_k + \int \Phi^{1-\nu}\zeta_k{\Gamma_{\bf S}}(\Phi,{L_{\bf S}}\Phi)+o(1).
\end{equation}
For $I_2$ at last, the chain rule simply implies that
\begin{equation}\label{I2}
    I_2=(\nu-1)\int\Phi^{1-\nu}\zeta_k{\Gamma_{\bf S}}(\Phi,{L_{\bf S}}\Phi).
\end{equation}
Now we turn to $II$ and apply the product rule.
$$
\frac2\nu II = \int\zeta_k {\Gamma_{\bf S}}({\Gamma_{\bf S}}(\Phi),\Phi^{1-\nu}) + \int {\Gamma_{\bf S}}(\Phi){\Gamma_{\bf S}}(\zeta_k,\Phi^{1-\nu})
= \int\zeta_k {\Gamma_{\bf S}}({\Gamma_{\bf S}}(\Phi),\Phi^{1-\nu}) + o(1).
$$
Thanks to the product rule again and integration by parts, it follows that
\begin{multline*}
\frac2\nu II = \int{\Gamma_{\bf S}}({\Gamma_{\bf S}}(\Phi),\zeta_k \Phi^{1-\nu}) + \int\Phi^{1-\nu}{\Gamma_{\bf S}}({\Gamma_{\bf S}}(\Phi),\zeta_k )=\\
\int{\Gamma_{\bf S}}({\Gamma_{\bf S}}(\Phi),\zeta_k \Phi^{1-\nu})-\int{\Gamma_{\bf S}}(\Phi){\Gamma_{\bf S}}(\Phi^{1-\nu},\zeta_k )-\int{\Gamma_{\bf S}}(\Phi)\Phi^{1-\nu}{L_{\bf S}}\zeta_k.
\end{multline*}
And so, since $\zeta_k\to1$ in $D({L_{\bf S}})$ and $\Phi,{\Gamma_{\bf S}}(\Phi)$ are bounded,
\begin{equation}\label{II}
\frac2\nu II=\int{\Gamma_{\bf S}}({\Gamma_{\bf S}}(\Phi),\zeta_k \Phi^{1-\nu}) + o(1) = -\int {L_{\bf S}}({\Gamma_{\bf S}}(\Phi))\zeta_k\Phi^{1-\nu}+o(1).    
\end{equation}
Plugging \eqref{I3}, \eqref{I1}, \eqref{I2}, \eqref{II} in \eqref{rhs}, we find
\begin{equation*}
\int\PAR{{\Gamma^{\bf S}_2}(\Phi)-\frac 1{\nu}({L_{\bf S}} \Phi)^2-\frac c{\nu}{\Gamma_{\bf S}}(\Phi)}\Phi^{1-{\nu}}\zeta_k = o(1).
\end{equation*}
Since $\Phi\in D(L_{\bf S})$, thanks to Lemma \ref{gammadeux}, we may pass to the limit as $k\to+\infty$ and deduce that
\begin{equation}
\label{eq:gamma2}
\int\PAR{{\Gamma^{\bf S}_2}(\Phi)-\frac 1{\nu}({L_{\bf S}} \Phi)^2-\frac c{\nu}{\Gamma_{\bf S}}(\Phi)}\Phi^{1-{\nu}} = 0.
\end{equation}
By the integrated curvature-dimension condition \eqref{icd}, we deduce that 
$$
\left(\frac 1n -\frac 1{\nu}\right)  \int({L_{\bf S}} \Phi)^2\Phi^{1-{\nu}} + \left(\rho-\frac c{\nu}\right)\int {\Gamma_{\bf S}}(\Phi)\Phi^{1-{\nu}}\le 0.
$$
Since $q<p$, we have $n< {\nu}$ and so, if  $\rho\ge\frac c{\nu}$ i.e.
$$
A\ge \frac{4({\nu}-1)}{{\nu}({\nu}-2)\rho},
$$
we deduce that ${L_{\bf S}} \Phi =0$. Integrating against $\Phi$, $\Phi$ is constant. Hence $v=1$, $I(A)=1$, and~\eqref{eq:sobno} holds for  $A=\frac{4({\nu}-1)}{{\nu}({\nu}-2)\rho}$.
Let $q\nearrow 2^*$. Then ${\nu}\searrow n$ and~the sharp inequality \eqref{355} follows.

It remains to study the case of equality. If $v\in H^1_0({\mu_{\bf S}})$ is an extremal function for~\eqref{355}, then repeating the above considerations, the function $f=v^{-\frac{p-2}2}$ satisfies
$$
\int \PAR{{\Gamma^{\bf S}_2}(f)-\rho{\Gamma_{\bf S}}(f)-\frac{1}{n}({L_{\bf S}} f)^2}f^{1-n}d{\mu_{\bf S}}= 0,
$$
In particular, if the parameters are such that inequality \eqref{opt2} is strict, it follows from \eqref{optr} that $f$ must be a function of $s$ only. If $\frac{n-d}{n-1}+\mathfrak B_{DGZ}= 0$, then 
the estimate~\eqref{687} in Lemma~\ref{lem-23} provides the following improvement of \eqref{optr}:
\begin{multline*}
\int \PAR{{\Gamma^{\bf S}_2}(f)-\rho{\Gamma_{\bf S}}(f)-\frac{1}{n}({L_{\bf S}} f)^2}f^{1-n}d{\mu_{\bf S}}
\geq \SBRA{\frac{n-d}{n-1}+\mathfrak B_{DGZ}}\int\Gamma^\theta(f)\phi^4f^{1-n}d{\mu_{\bf S}}\\
+\frac{n-d}{n-1}A\int \frac{\Gamma^\theta(f)^2}{f^2}\phi^4f^{1-n}d{\mu_{\bf S}}.
\end{multline*}
And so again, $f$ is a function of $s$ only, provided $n>d$.
Using this information in \eqref{554}, we deduce that if $n>d$,
$$
\partial_{ss}f +2\frac{\varphi'}{\varphi}\partial_sf=0,
$$
while for $n=d$ there must exist some function $R:\R\to\R$ s.t.
\begin{equation}\label{sphere}
\partial_{ss}f +2\frac{\varphi'}{\varphi}\partial_sf-\frac1{d-1}\Delta_\theta f=0\quad\text{and}\quad \partial_s f+\frac{\varphi'}{\varphi}f= R(s).
\end{equation}
In the former case, this means that $f(s)=\lambda+\gamma\tanh(\alpha s)$, for some constants $\lambda,\gamma\in\R$ such that $\lambda>\vert \gamma\vert$, since $f$ is bounded below by a positive constant. In the latter case, the second equation in \eqref{sphere} implies that $f$ can be written as  $f=\frac{f_1(\theta)}{\varphi(s)}+ f_2(s)$. Plugging this in the first equation implies that $f_1+\frac{\Delta_\theta f_1}{d-1}$ is constant i.e. $f_1=A_1+B_1\psi_2(\theta)$, where $A_1,B_1$ are constants and $\psi_2$ is any eigenfunction of $-\Delta_\theta$ associated to the eigenvalue $d-1$. This implies in turn that $f_2$ takes the form $f_2=-\frac{A_1}{\varphi(s)}+A_3+A_4\tanh(s)$. Summarizing, we have just proved that $f=\lambda+\gamma\varphi_{1,d}$ for some constants $\lambda,\gamma$ and some eigenfunction $\varphi_{1,d}$ of $-\Delta_{\Sp^d}$ associated to the eigenvalue $d$ (and written in cylindrical coordinates). Again, we must have $\lambda>\vert\gamma\vert \Vert\varphi_{1,d}\Vert_\infty$ since $f$ is bounded below by a positive constant.

Conversely, we need to check that $f^{-\frac{n-2}{2}}$ where $f(s)=\lambda+\gamma\tanh(\alpha s)$ with $\lambda>\vert\gamma\vert$ if $n>d$ (resp. $f=\lambda+\gamma\varphi_{1,d}$, $\lambda>\vert\gamma\vert \Vert\varphi_{1,d}\Vert_\infty$ if $n=d$) is indeed an extremal function for Sobolev's inequality. Multiplying $f$ by a constant if necessary, we may assume that $\int f^{-n}d\mu=1$, where  $\mu$ is the normalized measure on the CKN sphere. By direct computation, recalling that $\tanh(\alpha s)$ if $n>d$ (resp. $\varphi_{1,d}$ if $n=d$) is an eigenfunction for the operator $-{L_{\bf S}}$ associated to the eigenvalue $n\alpha^2$, we find that
$$
f{L_{\bf S}}f - \frac n2{\Gamma_{\bf S}}(f) = \frac{n\alpha^2}2(1-f^2)
$$
This implies in turn that  $v=f^{-\frac{n-2}{2}}$ satisfies $\int v^pd\mu=1$ and solves
$$
-\frac{4}{n(n-2)\alpha^2}{L_{\bf S}}v + v = v^{p-1}
$$
Multiplying by $v$ and integrating by parts, the result follows.
\end{eproofof}

\begin{eproofof}{Proposition \ref{poincare}}
As explained in the introduction, Poincaré's inequality (with constant $C=\frac{n-1}{\rho n}=\frac1{n\alpha^2}$) follows from Sobolev's inequality by linearization i.e. by applying~\eqref{355} with $v=1+\epsilon f$ and letting $\epsilon\to0$. Also, Poincaré's inequality (with the same constant $C$) is equivalent to the following integrated curvature-dimension condition
$$
\int \left(\Gamma_2^{\mathbf S}(f)- \rho\Gamma_{\mathbf S}(f)-\frac1n(L_{\mathbf S}f)^2\right)d\mu_{\mathbf S}\ge 0.
$$ 
Equality holds in Poincaré's inequality for some function $f$ if and only if equality holds in the above inequality. So, extremals are characterized exactly as in the case of Sobolev's inequality except in the case $\alpha^2=\frac{d-1}{n-1}$, in which we can no longer use \eqref{687} to deduce that $f$ is radial. Still, we deduce from \eqref{554} that
$$
\partial_{ss}f +2\frac{\varphi'}{\varphi}\partial_sf-\frac1{n-1}\Delta_\theta f=0\quad\text{and}\quad \partial_s f+\frac{\varphi'}{\varphi}f= R(s)
$$
The second equation in \eqref{sphere} implies that $f$ can be written as  $f=\frac{f_1(\theta)}{\varphi(s)}+ f_2(s)$. Plugging this in the first equation implies that $\alpha^2 f_1+\frac{\Delta_\theta f_1}{n-1}$ is constant i.e. $f_1=A_1+B_1\varphi_{1,d-1}(\omega)$, where $A_1,B_1$ are constants and $\varphi_{1,d-1}$ is any eigenfunction of $-\Delta_\theta$ associated to the eigenvalue $\alpha^2(n-1)=d-1$. This implies in turn that $f_2$ takes the form $f_2=-\frac{A_1}{\varphi(s)}+A_3+A_4\tanh(\alpha s)$. Summarizing, we have just proved that extremals of Poincaré's inequality take the form $f=\lambda+\gamma\tanh(\alpha s)+\nu\frac{\varphi_{1,d-1}(\omega)}{\cosh(\alpha s)}$ for some constants $\lambda,\gamma,\nu$ and some eigenfunction $\varphi_{1,d-1}$ of $-\Delta_{\Sp^{d-1}}$, as desired.
\end{eproofof}

\begin{erem}
Up to our knowledge, the CKN sphere is the first example where the optimal constants for both the Sobolev and the Poincaré inequalities are explicit functions of $(\rho,n)$ yet the usual curvature-dimension condition doesn't hold, although the integral version \eqref{cdrhoi} remains true. Beware though that the integrated curvature-dimension needed for (and equivalent to) Poincaré's inequality, i.e. inequality \eqref{cdrhoi} without the weight $f^{1-n}$, is in general much weaker, as evidenced by any space for which the Poincaré inequality holds but not the Sobolev inequality, such as, for instance, the Euclidean space equipped with the Gaussian measure.
\end{erem}

\section{The \boldmath{$n$}-conformal invariant}
\label{sec-6}
\subsection{The \boldmath{$n$}-conformal invariant on a weighted manifold}
We begin this section by proving Proposition \ref{prop1}, which constructs a one-parameter family of $n$-conformal invariants on any given weighted manifold, thereby generalizing the notion of scalar curvature to this setting.
 
\begin{eproofof}{Proposition \ref{prop1}}
We want to check that $S_\gamma(\mu,\Gamma)$ satisfies  condition~\eqref{178}. Let $c$ be a positive and smooth function on $M$,  $\tau=\log c$ and $\gamma\in\R$.  
We are looking for the expression of the two numbers $\theta_n(\gamma)$ and $\beta_n(\gamma)$ in the definition of $S_\gamma(\mu,\Gamma)$ which are such that 
$$
S_\gamma(c^{-n}\mu,c^2\Gamma)=c^2\SBRA{S_\gamma(\mu,\Gamma)+\frac{n-2}{2}\PAR{L\tau-\frac{n-2}{2}\Gamma(\tau)}}. 
$$
The measure $\mu$ is transformed into $\breve\mu=c^{-n}\mu$, and the carré du champ $\Gamma$ into $\breve\Gamma=c^{2}\Gamma$. 

From~\eqref{11},   $sc_{\mathfrak g}$ becomes 
$$
 {\breve {sc}}_{\mathfrak g}=c^2[sc_{\mathfrak g}+(d-1)(2\Delta_{\mathfrak g}\tau-(d-2)\Gamma(\tau))],
$$
the weight $W=-\log\frac{d\mu}{dV_{\mathfrak g}}$ becomes 
$$
\breve W=-\log\frac{d\breve\mu}{d\breve{V_{\mathfrak{g}}}}=-\log\frac{c^{-n}d\mu}{c^{-d}d{V_{\mathfrak{g}}}}=-\log \left(c^{d-n}\frac{d\mu}{d{V_{\mathfrak{g}}}}\right)=W+(n-d)\tau,
$$
and finally, from~\eqref{12}, $\Delta_{\mathfrak g}$ becomes 
$$
\breve\Delta_{\mathfrak g}=c^2[\Delta_{\mathfrak g}-(d-2)\Gamma(\tau,\cdot)].
$$

So, 
\begin{multline*}
S_\gamma(c^{-n}\mu,c^2\Gamma)=c^2\theta_n(\gamma)\Big[sc_{\mathfrak g}+(d-1)(2\Delta_{\mathfrak g}(\tau)-(d-2)\Gamma(\tau)) \\ -\gamma[\Delta_{\mathfrak g}(W+(n-d)\tau)-(d-2)\Gamma(\tau,W+(n-d)\tau)]+\beta_n(\gamma)\Gamma(W+(n-d)\tau)\Big]
\end{multline*}
that is 
\begin{multline*}
S_\gamma(c^{-n}\mu,c^2\Gamma)=c^2\theta_n(\gamma)\Big[sc_{\mathfrak g}+[2(d-1)-\gamma(n-d)]\Delta_{\mathfrak g}(\tau)\\
+[\beta_n(\gamma)(n-d)^2-(d-1)(d-2)+\gamma(d-2)(n-d)]\Gamma(\tau) \\ 
-\gamma\Delta_{\mathfrak g}(W) +[\gamma(d-2)+2\beta_n(\gamma)(n-d)]\Gamma(\tau,W)+\beta_n(\gamma)\Gamma(W)\Big].
\end{multline*}
It has to be equal to 
\begin{multline*}
c^2\Big[S_\gamma(\mu,\Gamma)+\frac{n-2}{2}\Big(\Delta_{\mathfrak g}(\tau)-\Gamma(W,\tau)-\frac{n-2}{2}\Gamma(\tau)\Big)\Big]=\\
c^2\Big[\theta_n(\gamma)[sc_{\mathfrak g}-\gamma\Delta_{\mathfrak g}(W)+\beta_n(\gamma)\Gamma(W)]+\frac{n-2}{2}\Big(\Delta_{\mathfrak g}(\tau)-\Gamma(W,\tau)-\frac{n-2}{2}\Gamma(\tau)\Big)\Big],
\end{multline*}
that is 
\begin{equation}
\label{61}
\left\{
\begin{array}{l}
\displaystyle\theta_n(\gamma)[2(d-1)-\gamma(n-d)]=\frac{n-2}{2} \\
\displaystyle\theta_n(\gamma)[\beta_n(\gamma)(n-d)^2-(d-1)(d-2)+\gamma(d-2)(n-d)]=-\frac{(n-2)^2}{4}\\
\displaystyle\theta_n(\gamma)[\gamma(d-2)+2\beta_n(\gamma)(n-d)]=-\frac{n-2}{2}
\end{array}
\right.
\end{equation}
which imples that 
\begin{equation*}
\left\{
\begin{array}{l}
\displaystyle\theta_n(\gamma)=\frac{n-2}{4(d-1)-2\gamma(n-d)} \\
\displaystyle\beta_n(\gamma)=\frac{\gamma(n-2d+2)-2(d-1)}{2(n-d)}.
\end{array}
\right.
\end{equation*}
Let us notice that the second equation in~\eqref{61} is automatically valid for this choice of parameters $\theta_n(\gamma)$ and $\beta_n(\gamma)$ and so we are done.
\end{eproofof}

\begin{erem}
As explained in the introduction, when $W=0$ the $d$-conformal invariant is, up to a multiplicative constant, the scalar curvature. In a weighted Riemannian manifold, the $n$-conformal invariant is given by~\eqref{181} and is a way to extend the definition of the scalar curvature in the weighted case.  	
\end{erem}
\subsection{The \boldmath{$n$}-conformal invariant for the CKN spaces}
In this section, we would like to prove that the three CKN spaces enjoy, for some $\gamma\in\R$, a {\it constant} $n$-conformal invariant. By construction, the three CKN models (Euclidean, spherical and hyperbolic) belong to the same $n$-conformal class. 
So, in virtue of Theorem~\ref{confinvckn} and Proposition~\ref{prop1}, it suffices to prove that there exits a unique $\gamma\in\R$ such that $S_\gamma=0$ for the Euclidean CKN space in order to prove Proposition~\ref{prop-5}.

\begin{eproofof}{Proposition~\ref{prop-5}}
Let $\gamma\in\R$. Then, 
$$
S_\gamma({\mu_{\bf E}},{\Gamma_{\bf E}})=\theta_n(\gamma)\big(sc_{{\mathfrak g_{\bf E}}}-\gamma\Delta_{{\mathfrak g_{\bf E}}}  W_{\bf E}+\beta_n(\gamma){\Gamma_{\bf E}}( W_{\bf E})\big).
$$
So, we need to find $\gamma$ such that 
$sc_{{\mathfrak g_{\bf E}}}-\gamma\Delta_{{\mathfrak g_{\bf E}}} W_{\bf E}+\beta_n(\gamma){\Gamma_{\bf E}}( W_{\bf E})=0$.

$\bullet$ Computation of the scalar curvature. From the identity~\eqref{11} with $\Gamma_{\bf E}=c_{\bf E}\Gamma$ and $\tau_{\bf E}=\log c_{\bf E}$, 
$$
sc_{{\mathfrak g_{\bf E}}}=|x|^{2(1-\alpha)}(0+(d-1)(2\Delta \tau_{\bf E}-(d-2)|\nabla \tau_{\bf E}|^2),
$$
hence, 
$$
sc_{{\mathfrak g_{\bf E}}}=|x|^{-2\alpha}(d-1)(d-2)(1-\alpha^2).
$$

$\bullet$ Computation of $\Delta_{{\mathfrak g_{\bf E}}} W_{\bf E}$. First, from the identity~\eqref{12}, 
$$
\Delta_{{\mathfrak g_{\bf E}}}  W_{\bf E}=|x|^{2(1-\alpha)\kappa}(\Delta  W_{\bf E}-(d-2)\nabla\tau_{\bf E}\cdot\nabla  W_{\bf E}),
$$
so 
$$
\Delta_{{\mathfrak g_{\bf E}}} W_{\bf E}=	|x|^{-2\alpha}(d-2)\alpha^2(n-d).
$$

$\bullet$ Computation of ${\Gamma_{\bf E}}(W_{\bf E})$. We have  
$$
{\Gamma_{\bf E}}( W_{\bf E})=|x|^{-2\alpha}\al^ha^2(n-d)^2.
$$
So, in the end, 
$$
S_\gamma(\mu,\Gamma)=\theta_n(\gamma)|x|^{-2\alpha}\PAR{(d-1)(d-2)(1-\alpha^2)+\gamma \alpha^2(n-d)(d-2)+\beta_n(\gamma)\alpha^2(n-d)^2},
$$
and we need to find $\gamma\in\R$ such that 
$$
(d-1)(d-2)(1-\alpha^2)+\gamma \alpha^2(n-d)(d-2)+\beta_n(\gamma)\alpha^2(n-d)^2=0.
$$
Since 
$$
\beta_n(\gamma)=\frac{\gamma(n-2d+2)-2(d-1)}{2(n-d)}
$$
we have 
\begin{equation*}
\gamma=2\frac{(d-1)(d-2+\alpha^2(2-n))}{\alpha^2(n-d)(2-n)}.
\end{equation*}
or by using the constant $\mathfrak B_{DGZ}$, 
\begin{equation*}
\gamma=\frac{2(d-1)}{\alpha^2(n-d)(2-n)}\mathfrak B_{DGZ}.
\end{equation*}
\end{eproofof}

\begin{erem}
It is interesting to notice that the $n$-conformal invariant for the CKN spaces does not depend on the sign of $\mathfrak B_{DGZ}$ or the Felli-Schneider region. 
\end{erem}

\begin{appendix}
\section{Appendix}

\subsection{Some Riemannian formulas}
\label{sec-a}
We recall here some general formulas on conformal transformations of a $d$-dimensional Riemannian manifold $(M,\mathfrak g)$. All formulas can be found for example in \cite[Sec.~6.9]{bgl-book}\footnote{and also here https://en.wikipedia.org/wiki/List\_of\_formulas\_in\_Riemannian\_geometry}. We transform the metric (with upper indices) $\mathfrak g^{ij}$ into the conformal metric ${\mathfrak h^{ij}}=c^2{\mathfrak  g^{ij}}$, where $c$ is any positive and smooth function. We let $\tau=\log c$. Then, 
\begin{itemize}
	\item The carré du champ operator is given by
	$$
	\Gamma^{ {\mathfrak  h}}=c^2\Gamma^{\mathfrak  g}.
	$$
	\item The Laplace-Beltrami operator is given by
	\begin{equation}
	\label{12}
	\Delta_{{\mathfrak  h}}=c^2(\Delta_{\mathfrak g}-(d-2)\Gamma^{\mathfrak  g}(\tau,\cdot)).
	\end{equation}
		\item For any smooth function $\psi$, the Hessian of $\psi$ with respect to the metric ${\mathfrak  h}$, denoted ${\nabla\nabla}^{\mathfrak h}  \psi$ is given by
	\begin{equation}
	\label{130}
	{\nabla\nabla}^{\mathfrak h}  \psi=\nabla\nabla^{\mathfrak g} \psi+2\nabla^{\mathfrak g}  \psi \odot_{\mathfrak g}\nabla^{\mathfrak g} \tau-\Gamma^{\mathfrak g}(\psi,\tau)
	{\mathfrak g},
	\end{equation}
	Here and below, $\nabla\nabla^{\mathfrak g} \psi$ is the Hessian of $\psi$ with respect to $\mathfrak g$ 
	 and $\nabla^{\mathfrak g}  \psi\odot_{\mathfrak g}\nabla^{\mathfrak g} \tau$ is the symmetric tensor product, that is for any functions $f,g$, 
	$$
	(\nabla^{\mathfrak g}  \psi\odot_{\mathfrak g}\nabla^{\mathfrak g} \tau)(\nabla^{\mathfrak g} f,\nabla^{\mathfrak g} g)=\frac{1}{2}\big[\Gamma^{\mathfrak g}(f, \psi)\Gamma^{\mathfrak g}(g,\tau )+\Gamma^{\mathfrak g}(f,\tau)\Gamma^{\mathfrak g}(g, \psi)\big].
	$$
	In particular, one can deduce the Hilbert-Schmidt norm of ${\nabla\nabla}^{\mathfrak h}  \psi$ with respect to the new metric $\mathfrak h$: 
	\begin{equation}
	\label{131}
	||{\nabla\nabla}^{\mathfrak h}  \psi||^2=c^4\Big[||\nabla\nabla^{\mathfrak g} \psi||^2+2\Gamma^{\mathfrak g}( \tau,\Gamma^{\mathfrak  g}(\psi))+2\Gamma^{\mathfrak  g}(\psi)\Gamma^{\mathfrak  g}(\tau)+
	(d-2)\Gamma^{\mathfrak  g}(\psi,\tau)^2-2(\Delta_{\mathfrak g} \psi)\Gamma^{\mathfrak  g}(\psi,\tau)\Big].
	\end{equation}
	  
	\item The Ricci curvature reads
	\begin{equation}
	\label{120}
	Ric_{{\mathfrak  h}}=Ric_{\mathfrak  g}+(\Delta_{\mathfrak  g}\tau){\mathfrak  g}+(d-2)(\nabla\nabla^{\mathfrak g}\tau+\nabla\tau\odot_{\mathfrak g}\nabla\tau-\Gamma^{\mathfrak g}(\tau){\mathfrak  g})
	\end{equation}
	
	\item At last, the scalar curvature is given by
	\begin{equation}
	\label{11}
	sc_{\mathfrak h}=c^2\big[sc_{\mathfrak g}+(d-1)\left(2\Delta_{\mathfrak g}\tau-(d-2)\Gamma^{\mathfrak g}(\tau)\right)\big].
	\end{equation}
\end{itemize}

\subsection{Integration by parts and elliptic theory on the CKN spherical space}\label{sec-a2}
%

Let $H^1_0({\mu_{\bf S}})$ denote the closure of $\mathcal C^\infty_c(\R^d\setminus\{0\})$ with respect to the norm
$$
\Vert u\Vert_{H^1_0({\mu_{\bf S}})}^2 = \int ({\Gamma_{\bf S}}(u)+u^2)\;d{\mu_{\bf S}}
$$
Let $u\in L^2({\mu_{\bf S}})$. Then,  $\vert x\vert^{1-\alpha}\frac{1+\vert x\vert^{2\alpha}}2\nabla u$ and ${L_{\bf S}} u$ are well-defined distributions on $\R^d\setminus\{0\}$ and we may ask whether they are actual functions in $L^2({\mu_{\bf S}})$, that is, we may consider ${L_{\bf S}}$ as an unbounded operator in $L^2({\mu_{\bf S}})$ with domain 
$$D({L_{\bf S}})=\{u\in H^1_0({\mu_{\bf S}})\;:\; {L_{\bf S}} u\in L^2({\mu_{\bf S}})\},$$
equipped with the norm 
$$\Vert u\Vert^2=\Vert u\Vert_{L^2({\mu_{\bf S}})}^2+\Vert {L_{\bf S}}u\Vert_{L^2({\mu_{\bf S}})}^2=\int \left(u^2+({L_{\bf S}}u)^2\right)\;d{\mu_{\bf S}}.$$
Since $\int {\Gamma_{\bf S}}(u)\;d{\mu_{\bf S}}=-\int uL_{\bf S}u\;d{\mu_{\bf S}}$ for $u\in\mathcal C^\infty_c(\R^d\setminus\{0\})$, it easily follows that $L_{\bf S}$ is a closed operator. In addition, the integration by parts formula holds on its domain:

\begin{elem}Let $u,v\in D({L_{\bf S}})$. Then,
\begin{equation}\label{IPP}
\int -{L_{\bf S}}u\, v\,d{\mu_{\bf S}} = \int {\Gamma_{\bf S}}(u,v)d{\mu_{\bf S}}
\end{equation}
\end{elem}
\begin{eproof}
Assume first that $u,v\in \mathcal C^\infty_c(\R^d\setminus\{0\})$. Then,~\eqref{IPP} follows by standard integration by parts. Next, if $u\in D({L_{\bf S}})$ and $v\in \mathcal C^\infty_c(\R^d\setminus\{0\})$, take $u_n\in \mathcal C^\infty_c(\R^d\setminus\{0\})$ s.t. $u_n\to u$ in $H^1_0(\mu)$. Using successively the definition of distributional derivatives, the convergence $u_n\to u$ in $L^2({\mu_{\bf S}})$, standard integration by parts and the convergence $u_n\to u$ in $H^1_0(\mu)$, we find
$$
\int -{L_{\bf S}} u\, v\,d{\mu_{\bf S}} = \int  u\, (-{L_{\bf S}}v)\,d{\mu_{\bf S}} = \lim_{n\to+\infty} \int u_n \, (-{L_{\bf S}}v)\,d{\mu_{\bf S}} = \lim_{n\to+\infty} \int {\Gamma_{\bf S}}(u_n,v)d{\mu_{\bf S}} = \int {\Gamma_{\bf S}}(u,v)d{\mu_{\bf S}}
$$
Finally, if  $u\in D({L_{\bf S}})$ and  if $v\in D({L_{\bf S}})$, take $v_n\in \mathcal C^\infty_c(\R^d\setminus\{0\})$ s.t. $v_n\to v$ in $H^1_0(\mu)$. Then, according to what we just proved,
$$
\int -{L_{\bf S}} u\, v\,d{\mu_{\bf S}} = \lim_{n\to+\infty} \int -{L_{\bf S}} u\, v_n\,d{\mu_{\bf S}} = \lim_{n\to+\infty} \int {\Gamma_{\bf S}}(u,v_n)d{\mu_{\bf S}} = \int {\Gamma_{\bf S}}(u,v)d{\mu_{\bf S}}
$$
\end{eproof}
The following approximation lemma will be useful to integrate by parts in more delicate settings than the above lemma. 
\begin{elem}\label{kesa}Assume $n>4$.
Let $u\in D(L_{\bf S})$ be such that $u,\Gamma_{\bf S}(u)$ are bounded. Then, there exists $u_k\in C^\infty_c(\R^d\setminus\{0\})$ such that $u_k\to u$ in $D_{L_{\bf S}}$.
\end{elem}
\begin{erem}The assumptions $n>4$ and $u,\Gamma_{\bf S}(u)$ bounded can be removed and replaced by $n\ge3$, but the proof is more involved (see \cite[Thm.~3.12]{Ketterer}).
\end{erem}
\begin{eproof} 
In Section~\ref{sec-3}, the model given in~\eqref{456} is written with the variables $(s,\theta)\in\R\times \Sp^{d-1}$. Choosing now $s=\frac{1}{\alpha}\mathrm{Argch}\PAR{1/\sin(t)}$ the model becomes with the variables $(t,\theta)\in(0,\pi)\times\Sp^{d-1}$,
\begin{equation}\label{sun}
    {L_{\bf S}}(f)=\alpha^2(\partial_{tt}f+c'(t)\partial_tf)+\frac{1}{\sin^2(t)}\Delta_\theta f,
\end{equation}

for any  smooth function $f$ defined in $(0,\pi)\times\Sp^{d-1}$, where $c(t)=(n-1)\log(\sin(t))$. The  carré du champ operator becomes  
\begin{equation}\label{sdeux}
{\Gamma_{\bf S}}(f)=\alpha^2(\partial_t f)^2+\frac{1}{\sin^2(t)}\Gamma^\theta(f)
\end{equation}
and invariant measure  
\begin{equation}\label{strois}
    d\mu_{\bf S}(t,\theta)=\frac{1}{Z}e^{c(t)}dtdV_{\bf S^{d-1}}=\frac{1}{Z}(\sin t)^{n-1}dt dV_{\bf S^{d-1}},
\end{equation}
where $Z$ is a normalization constant and $dV_{\bf S^{d-1}}$ is the volume in $\Sp^{d-1}$. 

Let now $\zeta_k\in C^\infty_c(0,\pi)$ denote a standard cut-off function such that $0\le \zeta_k\le 1$, $\zeta_k=0$ in $(0,1/k)\cup(\pi-1/k,\pi)$, $\zeta_k=1$ in $(2/k,\pi-2/k)$ and $|\zeta_k'|\le 2k$, $|\zeta_k''|\le 2k^2$. Setting 
$
u_k=u\zeta_k,    
$
we find
$$
\int (L_{\bf S}(u-u_k))^2 d\mu_{\bf S}\le 2 \int (L_{\bf S}u)^2(\zeta_k-1)^2d\mu_{\bf S}+ 4\int\Gamma_{\bf S}(u,\zeta_k)^2d\mu_{\bf S}+2\int u^2(L_S\zeta_k)^2d\mu_{\bf S}
=:I_1+I_2+I_3$$
By dominated convergence, $I_1\to0$ as $k\to+\infty$. For $I_2$, thanks to \eqref{sdeux} and \eqref{strois},
$$
I_2 \le C\Vert\Gamma_{\bf S}(u)\Vert_\infty^2k^{2-n}\to0
$$
Similarily, thanks to \eqref{sun} and \eqref{strois},
$$
I_3 \le C\Vert u\Vert_\infty^2k^{4-n}\to0.
$$
\end{eproof}

Our next tool is the following version of the Rellich-Kondrachov compactness theorem.
\begin{elem}\label{rellich}Let $(M,g,\mu)$ be a smooth connected weighted $d$-dimensional Riemannian manifold s.t. $d\ge3$, $\mu(M)<+\infty$ and Sobolev's inequality holds i.e. there exist constants $A,B\ge0$, $p\in\left[2,\frac{2d}{d-2}\right]$ such that for every $v\in \mathcal C^\infty_c(M)$,
\begin{equation}\label{non tight sobolev}
\left(\int \vert v\vert^pd\mu\right)^{\frac2p}\le A \int\Gamma(v)d\mu+B\int v^2d\mu
\end{equation}
Let $H^1_0(\mu)$ be the closure of $\mathcal A_0=\mathcal C^\infty_c(M)$ for the norm $\Vert u\Vert_{H^1_0(\mu)}^2 = \int (u^2+\Gamma(u))d\mu$ and let $q\in[1,p)$. Then, the embedding $H^1_0(\mu)\hookrightarrow L^q(\mu)$ is compact.
\end{elem}

\begin{eproof}Cover $M$ by a countable increasing family of open sets $(\Omega_k)_{k\in\N}$ with compact closure and for each $k\in\N$, let $\eta_k\in \mathcal C^\infty_c(\Omega_{k+1})$ be such that $\eta_k=1$ in $\Omega_k$. Let $(u_m)$ be a bounded sequence in $H^1_0(\mu)$. Since $d\mu=e^{-W}dV_g$ and $W,g$ are smooth, the $H^1_0(\mu)$ and the standard $H^1_0$ norm are equivalent for functions compactly supported in a fixed $\Omega_k$. By the classical Rellich-Kondrachov theorem, we deduce that for fixed $k$, the sequence $(u_m\eta_k)_m$ is compact in $L^r(\Omega_{k+1},d\mu)$ for $r\in(q,p)$. Since $(u_m)$ is bounded in the Hilbert space $H^1_0(\mu)$, by the Banach-Alaoglu theorem, $(u_m)$ is also compact in $H^1_0(\mu)$ for the weak topology. By a standard diagonal argument, a subsequence of $(u_m)$ (denoted the same) converges weakly in $H^1_0(\mu)$ to some function $u\in H^1_0(\mu)$ such that $(u_m\eta_k)_m$ converges to $u\eta_k$ in $L^r(\mu)$. 
Now, using H\"older's and Sobolev's inequality we find
\begin{align*}
\Vert u_m -u\Vert_{L^q(\mu)} &\le \Vert (u_m -u)\eta_k\Vert_{L^q(\mu)}+\Vert (u_m -u)(1-\eta_k)\Vert_{L^q(\mu)} \\
&\le 
 \Vert (u_m -u)\eta_k\Vert_{L^r(\mu)}\mu(M)^{\frac1q-\frac1r}+\Vert u_m -u\Vert_{L^p(\mu)}\mu(M\setminus\Omega_k)^{\frac1q-\frac1p}\\
 &\le C\left( \Vert (u_m -u)\eta_k\Vert_{L^r(\mu)}+\mu(M\setminus\Omega_k)\right)^{\frac1q-\frac1p})
\end{align*}
Hence,
$$
\limsup_{m\to+\infty}\Vert u_m -u\Vert_{L^q(\mu)}\le C\mu(M\setminus\Omega_k)^{\frac1q-\frac1p}
$$
Letting $k\to+\infty$, the claim follows.
\end{eproof}
As an immediate consequence of the above lemma (and a proof by contradiction), we have
\begin{ecor}\label{cor:poincare}
Make the same assumptions as in Lemma \ref{rellich}. Assume in addition that constants belong to $H^1_0(\mu)$. Then, Poincaré's inequality holds i.e. there exists a constant $C_P>0$ such that 
$$
\int v^2d\mu - \left(\int v d\mu\right)^2 \le C_P\int\Gamma(v)d\mu \qquad\text{for $v\in\mathcal A_0=\mathcal C^\infty_c(M)$}
$$
\end{ecor}
Finally, we state and prove elliptic regularity estimates, which are useful to justify integrations by parts in our proof of Sobolev's inequality.
\begin{elem}[General elliptic estimates]
\label{lem-975}
Assume that $(a,b)\in\Theta_{FS}$ (defined in~\eqref{865}). Let also  $h:(0,1]\times\Sp^{d-1}\mapsto\R$ be a smooth and bounded function satisfying $\int_{\Sp^{d-1}} hdV_{\bf S^{d-1}}=0$ and solving the equation 
$$
\partial_{tt}h-\frac{(n-1)(n-3)}{4}\frac{h}{t^2}+\frac{\Delta_\theta h}{\alpha^2 \sin^2(t)}=R,
$$
where $R$ is a smooth and bounded function on $(0,1]\times\Sp^{d-1}$. We assume also that, uniformly on $(0,1]\times\Sp^{d-1}$
$$
|h(t,\theta)|\leq Ct^{\frac{n-1}{2}}\quad {\rm and}\quad  |R(t,\theta)|\leq Ct^{\frac{n-1}{2}},
$$
for some constant $C>0$. Then, there exists a constant $C'>0$ s.t. uniformly on $(0,1]\times\Sp^{d-1}$,
$$
|h(t,\theta)|\leq C't^{\frac{n+1}{2}}, \quad \Gamma^\theta(h)(t,\theta)=|\nabla_\theta h|^2(t,\theta)\leq C't^{n+1}\quad {\rm and}\quad  
|\partial_th(t,\theta)|\leq C't^{\frac{n-1}{2}}.
$$
\end{elem}

\begin{eproof}
Let $(P_k)_{k\geq 0}$ be the orthonormal basis of eigenvectors of the operator $-\Delta_\theta$ on $\Sp^{d-1}$ associated to the increasing sequence of eigenvalues $(\lambda_k)_{k\geq 0}$ (recall that $\lambda_k\geq \lambda_1=d-1$, for $k\geq1$ and $\lambda_0=0$).

We decompose  $h$ in the basis $(P_k)_{k\geq 0}$,  
$$
h(t,\theta)=\sum_{k=1}^\infty h_k(t)P_k(\theta), \quad (t,\theta)\in(0,\pi)\times\Sp^{d-1},
$$
where $h_k(t)=\int_{\Sp^{d-1}}P_k (\theta)h(t,\theta)dV_{\bf S^{d-1}}$ (noting that  $\int hP_0dV_{\bf S^{d-1}}=\int hdV_{\bf S^{d-1}}=0$, whence $h_0=0$).  For each $k\geq 1$, $h_k$ satisfies 
$$
h_k''-\frac{(n-1)(n-3)}{4}\frac{h_k}{t^2}-\frac{\lambda_k}{\alpha^2\sin^2(t)}h_k=R_{k},
$$
where $R_{k}=\int_{\Sp^{d-1}} RP_kdV_{\bf S^{d-1}}$, which  satisfies again $|R_k(t,\theta)|\leq C_kt^{\frac{n-1}{2}}$.  The equation can be replaced by the following one
$$
h_k''-a\frac{h_k}{t^2}=R_{2,k},
$$
where $a=\frac{(n-1)(n-3)}{4}+\frac{\lambda_k}{\alpha^2}>0$ and $R_{2,k}$ satisfies the same estimate as $R_k$.

We are now able to solve the ODE. The method of variation of constants gives the explicit solution:
$$
h_k(t)=At^{\gamma_+}+Bt^{\gamma_-}+\frac{t^{\gamma_+}}{\gamma_+-\gamma_-}\int_1^t R_{2,k}(y)y^{1-\gamma_+}dy-\frac{t^{\gamma_-}}{\gamma_+-\gamma_-}\int_0^t R_{2,k}(y)y^{1-\gamma_-}dy,
$$
where $A$, $B$ are constants and
$$
\gamma_\pm=\frac{1\pm\sqrt{1+4a}}{2}=\frac{1}{2}\pm\frac{1}{2}\sqrt{(n-2)^2+4\frac{\lambda_k}{\alpha^2}},
$$
so that $\gamma_+>0$, $\gamma_-<0$.  
Then, by the estimate satisfied by $R_{2,k}$ we have near 0,
$$
\Big|\frac{t^{\gamma_+}}{\gamma_+-\gamma_-}\int_1^t R_{2,k}(y)y^{1-\gamma_+}dy\Big|\leq Ct^{\gamma_+}\Big|\int_1^t y^{\frac{n+1}{2}-\gamma_+}dy\Big|\leq  C\PAR{t^{\gamma_+}+t^{\frac{n+3}{2}}}
$$
and 
$$
\Big|\frac{t^{\gamma_-}}{\gamma_+-\gamma_-}\int_0^t R_{2,k}(y)y^{1-\gamma_-}dy\Big|\leq Ct^{\gamma_-}\int_0^t y^{\frac{n+1}{2}-\gamma_-}dy\leq Ct^{\frac{n+3}{2}}
$$
Since $h_k$ is a bounded function, we deduce that $B=0$ and
$$
|h_k(t)|\leq C(t^{\gamma_+}+t^{\frac{n+3}{2}}). 
$$
We claim that for $k\geq 1$,  $\gamma_+\geq \frac{n+1}{2}$, whence $|h_k(t)|\leq C t^{\frac{n+1}{2}}$. Indeed, by definition of $\gamma_+$, one can check that the inequality $\gamma_+\geq (n+1)/2$ is equivalent to  $\frac{\lambda_k}{n-1}\geq \alpha^2$. Since $\lambda_k\geq d-1$ for any $k\ge 1$ and since by Lemma~\ref{lem-15}, $(a,b)\in\Theta_{FS}$ if and only if $\alpha^2\leq\frac{d-1}{n-1}$, we indeed have $\gamma_+\geq (n+1)/2$. 

Next, we prove that the estimate remains valid for the function $h$. Define 
$$
H_K=\sum_{k=1}^K h_kP_k. 
$$
so that $h=\lim_{K\rightarrow\infty}H_K$ pointwise. From the previous computations, we know that, uniformly in $(0,1)\times\Sp^{d-1}$, 
$$
|H_K(t,\theta)|\leq C_K t^{\frac{n+1}{2}}
$$
%
We prove now that the inequality is uniform in the parameter $K$. Assume this is not the case i.e.
\begin{equation}
\label{775}  
\sup_{K\geq 1}||H_K||_{(\frac{n+1}{2})}=+\infty,
\end{equation}
where
$$
||H_K||_{(\frac{n+1}{2})}=\sup_{t\in(0,1], \theta\in\Sp^{d-1}}\frac{|H_K(t,\theta)|}{t^{\frac{n+1}{2}}}.
$$ 
There exist a sequence $((t_K,\theta_K))_{K\geq 1}$ in $(0,1]\times\Sp^{d-1}$,  such that 
\begin{equation}
\label{222}
\lim_{K\rightarrow\infty} \frac{|H_K(t_K,\theta_K)|}{t_K^\frac{n+1}{2}}=\infty. 
\end{equation}
By compactness, one can assume that $(t_K)$ (resp. $(\theta_K)$) converges to $t_\infty\in[0,1]$ (resp. $\theta_\infty\in\Sp^{d-1}$). There are two cases, either $t_\infty>0$ or $t_\infty=0$. The first case is not possible. Indeed, $h$ is bounded by assumption, whence $H_K$ is bounded by a constant independent of $K$ and so ~\eqref{222} contradicts $t_\infty>0$. 
The remaining case $t_\infty=0$ is more tricky. Let 
$$
G_K(z,\theta)=\frac{H_K(zt_K,\theta)}{t_K^{\frac{n+1}{2}}||H_K||_{(\frac{n+1}{2})}}
$$
for any $z\in(0,1/t_K]$. From the equation satisfied by each $h_k$, we have 
$$
\partial_{zz}G_K(z,\theta)-\frac{(n-1)(n-3)}{4}\frac{G_K(z,\theta)}{z^2}+t_ K^2\frac{\Delta_\theta G_K(z,\theta)}{\alpha^2\sin^2(z t_K)}=\frac{t_K^2}{t_K^{\frac{n+1}{2}}||H_K||_{(\frac{n+1}{2})}}\sum_{k=1}^KR_{k}(zt_k)P_k(\theta).
$$
By assumption on $R$, we have 
$$
\Big|\frac{t_K^2}{t_K^{\frac{n+1}{2}}||H_K||_{(\frac{n+1}{2})}}\sum_{k=1}^KR_{k}(zt_k)P_k(\theta)\Big|\leq C\frac{t_K}{||H_K||_{(\frac{n+1}{2})}}\stackrel{K\rightarrow\infty}{\longrightarrow} 0, 
$$
uniformly for $z$ in a compact subset of $\R^*$. By standard elliptic regularity, it follows that the sequence $(G_K)$ converges to $G$ solution on  $(0,\infty)\times \Sp^{d-1}$ of the PDE 
$$
\partial_{zz}G(z,\theta)-\frac{(n-1)(n-3)}{4}\frac{G(z,\theta)}{z^2}+\frac{\Delta_\theta G(z,\theta)}{\alpha^2 z^2}=0.
$$
Now, using the same argument we have   
$$
G=\sum_{k=0}^\infty G_kP_k.
$$
where again $G_k(z)=\int_{\Sp^{d-1}} G(z,\theta)P_k(\theta)dV_{\bf S^{d-1}}$. 
Then for each $k\geq 0$,  we have $G_k(t)=At^{\gamma_+}+Bt^{\gamma_-}$
where $\gamma_\pm$ are the same constants as before. But, the function $G_k$, defined on $(0,\infty)$, is bounded. This implies that  $A=B=0$ and then $G=0$. But by its definition, we know that  
$$
G(1,\theta_\infty)=1, 
$$
which gives a contradiction: the hypothesis~\eqref{775} is not valid. We conclude that  uniformly in $(0,1]\times \Sp^{d-1}$, we have  
$$
|h(t,\theta)|\leq Ct^{\frac{n+1}{2}}.
$$
It remains to prove the gradient estimates. Fix $t>0$ and for $z\in(1/4,2)$, $\theta\in \Sp^{d-1}$, let this time $G(z,\theta)=h(tz,\theta)$ and $S(z,\theta)=R(tz,\theta)$ so that
$$
\partial_{zz}G-\frac{(n-1)(n-3)}{4}\frac{G}{z^2}+\frac{t^2}{\alpha^2 \sin^2(tz)}\Delta_\theta G=t^2S\quad\text{in $(1/4,2)\times \Sp^{d-1}$}
$$
Note that the coefficients of the elliptic operator on the left-hand side are bounded in $C^2$-norm by a constant independent of $t$ so that, by standard elliptic regularity,
$$
\vert \partial_z G\vert+ \vert \nabla_\theta G\vert \le C\left(t^2\Vert S\Vert_{L^\infty((1/4,2)\times \Sp^{d-1})}+ \Vert G\Vert_{L^\infty((1/4,2)\times \Sp^{d-1})}  \right)\le C t^{\frac{n+1}2}\quad\text{in $(1/2,3/2)\times \Sp^{d-1}$,}
$$
for some constant $C>0$ independent of $t$. The desired estimates on $h$ follow by applying the above estimate at $z=1$.
\end{eproof}

\begin{elem}
\label{lem-985}
Whenever $(a,b)\in\Theta_{FS}$ (defined in~\eqref{865}), 
the solution $v$ of the equation~\eqref{10} has a bounded carré du champ operator 
$$
||{\Gamma_{\bf S}}(v)||_\infty<+\infty.
$$
\end{elem}

\begin{eproof}
We use the chart and notation introduced in \eqref{sun}. We have to prove that $v$, solution of~\eqref{10} has a bounded carré du champ, that is  $||{\Gamma_{\bf S}}(v)||_\infty<+\infty$.
Letting $c(t)=(n-1)\log(\sin(t))$ and $h=e^{\frac{c}{2}}v$, equation~\eqref{10} becomes 
\begin{equation}
\label{354}
\partial_{tt}h-\PAR{\frac{2c''+c'^2}{4}}h+\frac{1}{\alpha^2\sin^2(t)}\Delta_\theta h=R,
\end{equation}
where 
$$
R=\frac{1}{\alpha^2 A}\PAR{h-e^{\frac{c(2-q)}{2}}h^{q-1}}.
$$
This transformation allows us to deal with a simpler PDE. We know that $v=e^{-\frac{c}{2}}h$ is bounded and positive. So, for some constant $C$  (the value of which is allowed to change from line to line), 
$$
0\leq h\leq Ce^{\frac{c}{2}}=C\sin(t)^\frac{n-1}{2}.
$$
Thus $|R|\leq C \sin(t)^\frac{n-1}{2}$.    And then, from the definition of $h$, the following inequality, 
\begin{equation}
\label{219}
\Big(\partial_t h-\frac{c'}{2}h\Big)^2+\frac{\Gamma^\theta(h)}{\alpha^2\sin^2(t)}\leq C \sin(t)^{n-1}
\end{equation}
is equivalent to  $||{\Gamma_{\bf S}}(v)||_\infty<+\infty$. 

We know that $h$ is a smooth function on $(0,\pi)\times \Sp^{d-1}$. So, to prove the previous inequality, it is enough to work around $t=0$ and $t=\pi$. By symmetry, it suffices to treat the case $t=0$.

By definition of $c$, equation~\eqref{354} can be written as follow
$$
\partial_{tt}h-\frac{(n-1)(n-3)}{4}\frac{h}{\sin^2(t)}+\frac{1}{\alpha^2\sin^2(t)}\Delta_\theta h=R-\frac{(n-1)(n-3)}{4}h,
$$
or, since we are working around $t=0$, we have 
\begin{equation}
\label{547}
\partial_{tt}h-\frac{(n-1)(n-3)}{4}\frac{h}{t^2}+\frac{1}{\alpha^2\sin^2(t)}\Delta_\theta h=R_2,
\end{equation}
with 
$$
R_2=R-\frac{(n-1)(n-3)}{4}h+\frac{(n-1)(n-3)}{4}\PAR{\frac{1}{\sin^2(t)}-\frac{1}{t^2}}h,
$$
which satisfies again $|R_2|\leq C t^\frac{n-1}{2}$.

Let us write $h=h-h_0+h_0$ where $h_0=\int_{\Sp^{d-1}} hdV_{\bf S^{d-1}}$.  Then, since $h_0$ doesn't depend on $\theta$, 
$$
\Big(\partial_t h-\frac{c'}{2}h\Big)^2+\frac{\Gamma^\theta(h)}{\alpha^2\sin^2(t)}\leq 
2\Big(\partial_t (h-h_0)-\frac{c'}{2}(h-h_0)\Big)^2
+2\Big(\partial_t h_0-\frac{c'}{2}h_0\Big)^2+\frac{\Gamma^\theta(h-h_0)}{\alpha^2\sin^2(t)}.
$$
Then Lemma~\ref{lem-975} insures that $|h-h_0|\leq Ct^{\frac{n+1}{2}}$, $
|h'-h'_0|\leq Ct^{\frac{n-1}{2}}$ and $\Gamma^\theta(h-h_0)=\Gamma^\theta(h)\leq Ct^{n+1}$. Hence, 
$$
2\Big(\partial_t (h-h_0)-\frac{c'}{2}(h-h_0)\Big)^2+2\frac{\Gamma^\theta(h-h_0)}{\alpha^2\sin^2(t)}\leq C t^{n-1}.
$$
Now, using the same method as in the proof of Lemma~\ref{lem-975}, one can check that  $h_0(t)=A t^{\frac{n-1}{2}}+O(t^{\frac{n+1}{2}})$ and $h_0'=A t^{\frac{n-3}{2}}+O(t^{\frac{n-1}{2}})$. Thus, we have
$$
\partial_t h_0-\frac{c'}{2}h_0=O(t^{\frac{n-1}{2}})
$$
that is  
$$
\Big(\partial_t h_0-\frac{c'}{2}h_0\Big)^2\leq C t^{n-1}.
$$
Finally, inequality~\eqref{219} is satisfied, which concludes the proof. 
\end{eproof}

\begin{erem}
It is interesting to see that with this method, one can check that the function $v$ has a bounded carré du champ if and only of $(a,b)\in\Theta_{FS}$.
\end{erem}
We end this section with the following weaker estimate on higher derivatives of $v$.
\begin{elem} \label{gammadeux}
Assume that $n> 4$ and $(a,b)\in\Theta_{FS}$. Let $f\in D(L_{\mathbf S})$, $f\ge 0$. Then, in the variables $(t,\theta)\in(0,\pi)\times\Sp^{d-1}$ introduced in \eqref{sun},
$$
\int_0^\pi \left\vert\int_{\bf S^{d-1}}\Gamma_2^{\mathbf S}(f)dV_{\bf S^{d-1}}\right\vert \sin^{n-1}(t)dt<+\infty
$$
\end{elem}

\begin{eproof}
Let $(\zeta_k)_{k\in\N}\in C^\infty_c(0,\pi)$ denote a standard cut-off function such that $0\le \zeta_k\le 1$, $\zeta_k=0$ in $(0,1/k)\cup(\pi-1/k,\pi)$, $\zeta_k=1$ in $(2/k,\pi-2/k)$ and $|\zeta_k'|\le 2k$, $|\zeta_k''|\le 2k^2$ and $f_k=f\zeta_k$, so that $f_k\to f$ in $D(L_{\mathbf S})$. For $h\in C^\infty_c(\R^d\setminus\{0\})$ and $t\in(0,\pi)$, set $\gamma_2(h)(t)=\int_{\Sp^{d-1}}\Gamma_{\mathbf S}^2(h)dV_{\bf S^{d-1}}$. Since $(a,b)\in\Theta_{FS}$, $\gamma_2$ is a nonnegative quadratic form (see Proposition~\ref{crucial}) and so the Cauchy-Schwarz inequality holds:
$$
|\gamma_2(f_k)^{1/2}-\gamma_2(f_l)^{1/2}|^2\leq \gamma_2(f_k-f_l).
$$
Thus, letting $d\mu_t=\frac1Z\sin^{n-1}(t)dt$, 
$$
\int_0^\pi |\gamma_2(f_k)^{1/2}-\gamma_2(f_l)^{1/2}|^2d\mu_t\leq \int_0^\pi\gamma_2(f_k-f_l)d\mu_t=\int \Gamma_{\bf S}^2(f_k-f_l)d\mu_{\bf S}=\int (L_{\bf S}(f_k-f_l))^2d\mu_{\bf S}.
$$
Hence, $(\gamma_2(f_k)^{1/2})$ is a Cauchy sequence in $L^2(d\mu_x)$ and so $(\gamma_2(f_k))$ converges to some function $\gamma$ in $L^1(d\mu_t)$. In addition,
for fixed $t$, there exists $K= K_t$ such that for all $k\ge K$ and $\theta\in {\bf S}^{d-1}$, $f_k(t,\theta)=f(t,\theta)$, whence $\gamma_2(f_k)(t)\to\gamma_2(f)(t)$ for all $t\in(0,\pi)$. Hence, $\gamma=\gamma_2(f)$ and the lemma follows.
\end{eproof}

\subsection{List of constants and regions of parameters}
\label{sec-b}

We recall in this section the definition of the parameters and also some useful properties. Recall that $d\in\N$ is the topological dimension of the considered spaces, and that we assume that $d\geq 3$. Recall from the introduction the definition of the parameter range 
$$
\Theta=\{(a,b)\in\R^2,\,a\leq b< a+1,\,a<a_c\},
$$
where $a_c=(d-2)/2$. This is the set of parameters $(a,b)$  where the CKN inequality~\eqref{1} holds for all test functions $v\in \mathcal C^\infty_c(\R^d)$ which need not vanish near the origin  (recall that the limit case $b=a+1$ has been removed for simplicity). 
We also defined the number $\alpha=1+a-\frac{pb}{2}$, that is  
$$
\alpha=\frac{(a_c-a)(a+1-b)}{a_c-a+b}.
$$

Clearly $\alpha\geq 0$, for any $(a,b)\in\Theta$, including the limiting case $a=b=0$ for which $\alpha=1$. For any $(a,b)\in\Theta$, the exponent $p$ is given by 
$$
p=\frac{d}{a_c-a+b}<2^*=\frac{2d}{d-2}
$$
and 
$$
p=\frac{2n}{n-2},
$$
that is 
$$
n=d+\frac{d(b-a)}{1+a-b}=\frac{d}{1+a-b}.
$$
We always have $n\geq d$, and we call $n$ the intrinsic dimension of the considered model spaces. 
From a straightforward computation, we have   
\begin{equation}
\label{99}
d-n\alpha -pb=0
\end{equation}

The constant $\mathfrak B_{DGZ}=\alpha^2(2-n)+d-2$ which appears throughout the paper takes the following form with respect to $a$ and $b$: 
\[
\mathfrak B_{DGZ}=-2\frac{(a_c-a)^2(1+a-b)}{a_c+b-a}+2a_c
\]
Let us also recall the definition of the Felli-Schneider region: for $a\leq 0$, 
\begin{equation*}
b_{FS}(a)=\frac{d(a_c-a)}{2\sqrt{(a_c-a)^2+d-1}}+a-a_c,\quad\text{and }\Theta_{FS}=\{(a,b)\in\Theta,\, b\geq b_{FS}(a)\text{ if }a\leq 0\}
\end{equation*}
Let us prove Lemma \ref{lem-15}, which simplifies the expression of the Felli-Schneider region and shows its relation to our region $\Theta_{DGZ}=\left\{(a,b)\in \Theta,\, \mathfrak B_{DGZ}\geq0\right\}$.

\begin{eproofof}{Lemma \ref{lem-15}}
  Since 
  $$
  \mathfrak B_{DGZ}+\frac{n-d}{n-1}=(n-2)\PAR{ -\alpha^2+ \frac{d-1}{n-1}},
  $$
  we have 
  $$
  \left\{(a,b)\in \Theta,\, \mathfrak B_{DGZ}+\frac{n-d}{n-1}\geq0\right\}=\left\{(a,b)\in \Theta,\, \alpha^2\leq \frac{d-1}{n-1}\right\}.
  $$
  The fact that 
  $$
  \Theta_{FS}=\left\{(a,b)\in \Theta,\, \alpha^2\leq \frac{d-1}{n-1}\right\}
  $$
  is more delicate and is proved in~\cite[Sec.~3]{DolEstLos16}.

  Finally, the identity  $\Theta_{FS}=\{(a,b)\in\Theta,\,\alpha\in [0,1]\}$ is a little trickier.  Since $(a,b)\in\Theta$, then 
  \begin{equation}
    \label{415}
    0<\alpha=\frac{(a_c-a)(a+1-b)}{a_c-a+b}\leq 1-\frac{a}{a_c}.
  \end{equation}
  
  In the case where $a\geq0$ it follows that $0<\alpha\leq1$.   Assume now that $a<0$. By definition of $\alpha$, 
  $$
  b=\frac{(a_c-a)(a+1-\alpha)}{\alpha+a_c-a}.
  $$
  Then the inequality $b\geq\mathfrak{B}_{FS}$ is equivalent to 
  $$
  \frac{1+a_c}{\alpha+a_c-a}\geq \frac{d}{2\sqrt{(a_c-a)^2+d-1}}.
  $$
  Since $1+a_c=d/2$, the previous condition becomes
  $$
  a\geq \frac{\alpha^2+\alpha(d-2)-d+1}{2\alpha}=\frac{(\alpha-1+d)(\alpha-1)}{2\alpha}.
  $$
  If $\alpha>1$, this inequality implies $a>0$, contradicting our assumption.  We just proved that 
  $$
  \Theta_{FS}\subset \left\{(a,b)\in \Theta,\,\alpha\in(0,1]\right\}.
  $$
  Since $\alpha=0$ for $(a,b)=(0,1)$, $\alpha=1$ for $(a,b)=(0,0)$, since $\Theta$ is connected and since $\alpha$ depends continuously on $(a,b)$, we deduce that $\Theta_{FS}=\left\{(a,b)\in \Theta,\,\alpha\in(0,1]\right\}$. At last, if $\alpha=1$, then, from the previous displayed inequality, $a\geq 0$ and from~\eqref{415}, $a\leq0$. That is, $a=b=0$. 
\end{eproofof}

\begin{erem}\label{mufini}
The normalizing constant $Z$ defined in Theorem \ref{thm-1} is finite.
Indeed, using \eqref{99} and the change of variable $\vert x\vert =e^t$, we find
$$
Z=\int\! d{\mu_{\bf S}}=\int_{\R^d\setminus\{0\}}\!\!\! \left(\frac{2}{1+\vert x\vert^{2\alpha}}\right)^n\! \vert x\vert^{-bp}dx=\int_{\R^d\setminus\{0\}}\!\!\! \left(\frac{2}{\vert x\vert^{-\alpha}+\vert x\vert^{\alpha}}\right)^n \!\vert x\vert^{-d}dx
=\vert \Sp^{d-1}\vert\int_{\R} \!\!\cosh(\alpha t)^{-n}dt<\infty.
$$
\end{erem}

\end{appendix}

\bibliographystyle{alpha}
{\footnotesize{\bibliography{biblio}}}
\small
This work was supported by the French ANR-17-CE40-0030 EFI project. 

\medskip

L. D., I. G. {Institut Camille Jordan, Umr Cnrs 52065, Universit\'e Claude Bernard Lyon 1, 43 boulevard du 11 novembre 1918, F-69622 Villeurbanne cedex.
	\texttt{dupaigne, gentil@math.univ-lyon1.fr}
	
	\medskip
	
	S. Z.  ENS de Lyon,  CNRS, UMPA UMR 5669, F-69364 Lyon Cedex 07.
	\texttt{simon.zugmeyer@ens-lyon.fr}

\end{document}